\documentclass[12pt,oneside]{article}
\usepackage{amsmath,amssymb,amsfonts,amsthm}
\def\goth{\mathfrak}

\newcommand{\T}{{\cal T}}

\newcommand{\Real}{\mathbb R}

\newcommand{\To}{\longrightarrow}

\newcommand{\p}{\pi^{-1}(TM)}
\newcommand {\cp}{\mathfrak{X}(\pi (M))}
\newcommand {\cpp}{\mathfrak{X}(\T M)}

\def\pa{\partial}
\def\paa{\dot{\partial}}

\def\x{{\goth X}(\T M)}

\def\Section#1{\vspace{30truept}\addtocounter{section}{1}\setcounter{thm}{0}
\setcounter{equation}{0}{\noindent\Large\bf
    \arabic{section}.~~#1}\par \vspace{12pt}}

\newtheorem{thm}{Theorem}[section]
\newtheorem{cor}[thm]{Corollary}
\newtheorem{lem}[thm]{Lemma}
\newtheorem{prop}[thm]{Proposition}
\newtheorem{defn}[thm]{Definition}

\numberwithin{equation}{section}

\begin{document}
\title{\bf{A GLOBAL APPROACH TO THE THEORY OF CONNECTIONS IN FINSLER GEOMETRY } }
\author{{\bf Nabil L. Youssef$^{\dag}$, S. H. Abed$^{\dag}$ and A. Soleiman$^{\ddag}$}}
\date{}

\maketitle                     
\vspace{-1.15cm}
\begin{center}
{$^{\dag}$Department of Mathematics, Faculty of Science,\\ Cairo
University, Giza, Egypt}
\end{center}
\vspace{-0.8cm}
\begin{center}
nlyoussef2003@yahoo.fr,\ sabed52@yahoo.fr
\end{center}
\vspace{-0.7cm}
\begin{center}
and
\end{center}
\vspace{-0.7cm}
\begin{center}
{$^{\ddag}$Department of Mathematics, Faculty of Science,\\ Benha
University, Benha,
 Egypt}
\end{center}
\vspace{-0.8cm}
\begin{center}
soleiman@mailer.eun.eg
\end{center}
\smallskip

\vspace{1cm} \maketitle
\smallskip

\noindent{\bf Abstract.} Adopting the pullback approach to Finsler
geometry, the aim of the present paper is to provide intrinsic
(coordinate-free) proofs of the existence and uniqueness theorems
for the Chern (Rund) and Hashiguchi connections on a Finsler
manifold. To accomplish this, we introduce and investigate the
notions of semispray and nonlinear connection associated with a
given regular connection, in the pullback bundle. Moreover, it is
shown that for the the Chern (Rund) and Hashiguchi connections, the
associated semispray coincides with the canonical spray and the
associated nonlinear connection coincides with the Barthel
connection. Explicit intrinsic expressions relating these
connections and the Cartan connection are deduced.
\par
Although our investigation is entirely global, the local expressions
of the obtained results, when calculated, coincide with the existing
classical local results.
\par
We provide, for the sake of completeness and for comparison reasons,
two appendices, one of them presenting a global survey of canonical
linear connections in Finsler geometry and the other presenting a
local survey of our global approach.


\newpage


\vspace{30truept}\centerline{\Large\bf{Introduction}}\vspace{12pt}
\par
\par
The most well-known and widely used approaches to GLOBAL Finsler
geometry are the Klein-Grifone (KG-) approach (cf. \cite{r21},
\cite{r22}, \cite{r27}) and the pullback (PB-) approach (cf.
\cite{r58}, \cite{r61}, \cite{r74}, \cite{r44}). The universe of the
first approach is the tangent bundle of $\,\T M$ (i.e, $\pi_{\T
M}:T\T M\longrightarrow \T M$), whereas the universe of the second
is the pullback of the tangent bundle $TM$ by $\pi: \T
M\longrightarrow M$  (i.e., $P:\pi^{-1}(TM)\longrightarrow \T M$).
Each of the two approaches has its own geometry which differs
significantly from the geometry of the other (in spite of the
existence of some links between them).
\par
 The theory of connections is an important field of differential
 geometry. It was initially developed to solve pure geometrical
 problems. The most important linear connections in Finsler geometry were
 studied \textbf{locally} in \cite{r91}, \cite{r34}, \cite{r93},
 \cite{r95},...etc.
\par
In \cite{r92}, we have introduced and investigated new intrinsic
proofs of intrinsic versions of the existence and uniqueness
theorems for the Cartan and Berwald connections on a Finsler
manifold $(M,L)$. On the other hand, there are other connections of
particular important in Finsler geometry, namely Chern and
Hashiguchi connections.  To the best of our knowledge, there is no
proof of the existence and uniqueness theorems for the Chern and
Hashiguchi connections from a purely global perspective.
\par
The  main purpose of the present paper is to provide
\textbf{intrinsic} (coordinate-free) proofs of  the existence and
uniqueness theorems for the Chern and Hashiguchi connections within
the pullback formalism, making simultaneous use of some concepts and
results from the KG-approach. These proofs have the advantages of
being simple, systematic and parallel to and guided by the
Reimannian case.
\par
The  paper consists of three parts preceded by an introductory
section $(\S 1)$, which provides a brief account of the basic
definitions and concepts necessary for this work. For more details,
we refer to \cite{r44}, \cite{r61}, \cite{r21} and \cite{r22}.
\par
In the first part $(\S 2)$, we review the fundamental results
concerning the existence and uniqueness theorems for the Cartan  and
Berwald connections on a Finsler manifold $(M,L)$ \cite{r92}. From
these results, the relationships between the curvature tensors
associated with the Berwald connection $D^{\circ}$ and the Cartan
connection $\nabla$ are obtained.
\par
The second part $(\S 3)$ is devoted to an intrinsic proof of the
existence and uniqueness theorem of the Chern connection on a
Finsler manifold $(M,L)$ (Theorem \ref{th.r1}). For the Chern
connection, it is shown that the associated nonlinear connection
coincides with the Barthel connection (Theorem \ref{th.a}). This
establishes an important link between the PB-approach and the
KG-approach. Moreover, the relationship between this connection and
the Cartan connection is obtained (Theorem \ref{th.r5}).
\par
The third and last part $(\S 4)$ provides an intrinsic proof of the
existence and uniqueness theorem of the Hashiguchi connection on
$(M,L)$ (Theorem \ref{th.h2}).  The associated nonlinear connection
is shown to  coincide with the Barthel connection (Theorem
\ref{0th.h2}). As in the previous section, the relationship between
this connection and the Cartan connection is obtained (Theorem
\ref{th.t12}).
\par
We have to emphasize that without the insertion of the KG-approach,
we would have been unable to achieve these results. It should also
be pointed out that the present work is formulated in a prospective
modern coordinate-free form; the local expressions of the obtained
results, when calculated, coincide with the existing classical local
results.
\par
Finally, for the sake of completeness and for comparison reasons, we
provide two appendices, one of them presenting a global survey of
canonical linear connections in Finsler geometry and the other
presenting a local survey  of our global approach.
.


\Section{Notation and Preliminaries}

In this section, we give a brief account of the basic concepts
 of the pullback formalism necessary for this work. For more
 details, we refer to \cite{r58},\,\cite{r61},\,\cite{r74} and~\,\cite{r44}.
 We make the
assumption that the geometric objects we consider are of class
$C^{\infty}$.\\ The
following notation will be used throughout this paper:\\
 $M$: a paracompact real differentiable manifold of finite dimension $n$ and of
class $C^{\infty}$,\\
 $\mathfrak{F}(M)$: the $\Real$-algebra of differentiable functions
on $M$,\\
 $\mathfrak{X}(M)$: the $\mathfrak{F}(M)$-module of vector fields
on $M$,\\
$\pi_{M}:TM\longrightarrow M$: the tangent bundle of $M$,\\
$\pi: \T M\longrightarrow M$: the subbundle of nonzero vectors
tangent to $M$,\\
$V(TM)$: the vertical subbundle of the bundle $TTM$,\\
 $P:\pi^{-1}(TM)\longrightarrow \T M$ : the pullback of the
tangent bundle $TM$ by $\pi$,\\
 $\mathfrak{X}(\pi (M))$: the $\mathfrak{F}(\T M)$-module of
differentiable sections of  $\pi^{-1}(T M)$,\\
$ i_{X}$ : the interior product with respect to  $X
\in\mathfrak{X}(M)$,\\
$df$ : the exterior derivative  of $f$,\\
$ d_{L}:=[i_{L},d]$, $i_{L}$ being the interior derivative with
respect to a vector form $L$.

\par Elements  of  $\mathfrak{X}(\pi (M))$ will be called
$\pi$-vector fields and will be denoted by barred letters
$\overline{X} $. Tensor fields on $\pi^{-1}(TM)$ will be called
$\pi$-tensor fields. The fundamental $\pi$-vector field is the
$\pi$-vector field $\overline{\eta}$ defined by
$\overline{\eta}(u)=(u,u)$ for all $u\in \T M$.

We have the following short exact sequence of vector bundles,
relating the tangent bundle $T(\T M)$ and the pullback bundle
$\pi^{-1}(TM)$:\vspace{-0.1cm}
$$0\longrightarrow
 \pi^{-1}(TM)\stackrel{\gamma}\longrightarrow T(\T M)\stackrel{\rho}\longrightarrow
\pi^{-1}(TM)\longrightarrow 0 ,\vspace{-0.1cm}$$
 where the bundle morphisms $\rho$ and $\gamma$ are defined respectively by
$\rho := (\pi_{\T M},d\pi)$ and $\gamma (u,v):=j_{u}(v)$, where
$j_{u}$  is the natural isomorphism $j_{u}:T_{\pi_{M}(v)}M
\longrightarrow T_{u}(T_{\pi_{M}(v)}M)$. The vector $1$-form $J$ on
$TM$ defined by $J:=\gamma\circ\rho$ is called the natural almost
tangent structure of $T M$. The vertical vector field $\mathcal{C}$
on $TM$ defined by $\mathcal{C}:=\gamma\circ\overline{\eta} $ is
called the fundamental or the canonical (Liouville) vector field.

Let $D$ be  a linear connection (or simply a connection) on the
pullback bundle $\pi^{-1}(TM)$.
 We associate with
$D$ the map \vspace{-0.1cm}
$$K:T \T M\longrightarrow \pi^{-1}(TM):X\longmapsto D_X \overline{\eta}
,\vspace{-0.1cm}$$ called the connection (or the deflection) map of
$D$. A tangent vector $X\in T_u (\T M)$ is said to be horizontal if
$K(X)=0$ . The vector space $H_u (\T M)= \{ X \in T_u (\T M) :
K(X)=0 \}$ of the horizontal vectors
 at $u \in  \T M$ is called the horizontal space to $M$ at $u$  .
   The connection $D$ is said to be regular if
\begin{equation}\label{direct sum}
T_u (\T M)=V_u (\T M)\oplus H_u (\T M) \qquad \forall u\in \T M .
\end{equation}
\par If $M$ is endowed with a regular connection, then the vector bundle
   maps
\begin{eqnarray*}
 \gamma &:& \pi^{-1}(T M)  \To V(\T M), \\
   \rho |_{H(\T M)}&:&H(\T M) \To \pi^{-1}(TM), \\
   K |_{V(\T M)}&:&V(\T M) \To \pi^{-1}(T M)
\end{eqnarray*}
 are vector bundle isomorphisms.
   Let us denote
 $\beta:=(\rho |_{H(\T M)})^{-1}$,
then \vspace{-0.2cm}
   \begin{align}\label{fh1}
    \rho\circ\beta = id_{\pi^{-1} (TM)}, \quad  \quad
       \beta\circ\rho =\left\{
                                \begin{array}{ll}
                                          id_{H(\T M)} & {\,\, on\,\,   H(\T M)} \\
                                         0 & {\,\, on \,\,   V(\T M)}
                                       \end{array}
                                     \right.\vspace{-0.2cm}
\end{align}
The map $\beta$ will be called the horizontal map of the connection
$D$.
\par According to the direct sum decomposition (\ref{direct
sum}), a regular connection $D$ gives rise to a horizontal projector
$h_{D}$ and a vertical projector $v_{D}$, given by
\begin{equation}\label{proj.}
h_{D}=\beta\circ\rho ,  \ \ \ \ \ \ \ \ \ \ \
v_{D}=I-\beta\circ\rho,
\end{equation}
where $I$ is the identity endomorphism on $T(TM)$: $I=id_{T(TM)}$.
\par
 The (classical)  torsion tensor $\textbf{T}$  of the connection
$D$ is defined by
$$\textbf{T}(X,Y)=D_X \rho Y-D_Y\rho X -\rho [X,Y] \quad
\forall\,X,Y\in \mathfrak{X} (\T M).$$ The horizontal ((h)h-) and
mixed ((h)hv-) torsion tensors, denoted by $Q $ and $ T $
respectively, are defined by \vspace{-0.2cm}
$$Q (\overline{X},\overline{Y})=\textbf{T}(\beta \overline{X}\beta \overline{Y}),
\, \,\, T(\overline{X},\overline{Y})=\textbf{T}(\gamma
\overline{X},\beta \overline{Y}) \quad \forall \,
\overline{X},\overline{Y}\in\mathfrak{X} (\pi (M)).\vspace{-0.2cm}$$
If $M$ is endowed with a metric $g$ on $\p$, we write
\begin{equation}\label{tor.g}
    T(\overline{X},\overline{Y},\overline{Z}):
=g(T(\overline{X},\overline{Y}),\overline{Z}).
\end{equation}
\par
The (classical) curvature tensor  $\textbf{K}$ of the connection $D$
is defined by
 $$ \textbf{K}(X,Y)\rho Z=-D_X D_Y \rho Z+D_Y D_X \rho Z+D_{[X,Y]}\rho Z
  \quad \forall\, X,Y, Z \in \mathfrak{X} (\T M).$$
The horizontal (h-), mixed (hv-) and vertical (v-) curvature
tensors, denoted by $R$, $P$ and $S$ respectively, are defined by
$$R(\overline{X},\overline{Y})\overline{Z}=\textbf{K}(\beta
\overline{X}\beta \overline{Y})\overline{Z},\quad
P(\overline{X},\overline{Y})\overline{Z}=\textbf{K}(\beta
\overline{X},\gamma \overline{Y})\overline{Z},\quad
S(\overline{X},\overline{Y})\overline{Z}=\textbf{K}(\gamma
\overline{X},\gamma \overline{Y})\overline{Z}.$$ The contracted
curvature tensors, denoted by $\widehat{R}$, $\widehat{P}$ and
$\widehat{S}$ respectively, are also known as the
 (v)h-, (v)hv- and (v)v-torsion tensors and are defined by
$$\widehat{R}(\overline{X},\overline{Y})={R}(\overline{X},\overline{Y})\overline{\eta},\quad
\widehat{P}(\overline{X},\overline{Y})={P}(\overline{X},\overline{Y})\overline{\eta},\quad
\widehat{S}(\overline{X},\overline{Y})={S}(\overline{X},\overline{Y})\overline{\eta}.$$
If $M$ is endowed with a metric $g$ on $\p$, we write
\begin{equation}\label{cur.g}
    R(\overline{X},\overline{Y},\overline{Z}, \overline{W}):
=g(R(\overline{X},\overline{Y})\overline{Z}, \overline{W}),\,
\cdots, \, S(\overline{X},\overline{Y},\overline{Z}, \overline{W}):
=g(S(\overline{X},\overline{Y})\overline{Z}, \overline{W}).
\end{equation}

The following lemma is useful for subsequent use.\vspace{-0.2cm}
\begin{lem}\label{cyclic}For every linear connection $D$ on $\pi^{-1}(TM)$
with \textsc{(}classical\textsc{)} torsion $\textbf{T}$ and
\textsc{(}classical\textsc{)} curvature $\textbf{K}$, we have
\begin{description}
    \item[(a)]$\mathfrak{S}_{X,Y,Z}\{\textbf{K}(X,Y)\rho Z+D_{X}\textbf{T}(Y,Z)
    +\textbf{T}(X,[Y,Z])\}=0$,
    \item[(b)]$\mathfrak{S}_{X,Y,Z}\{D_{X}\textbf{K}(Y,Z)-
    \textbf{K}(X,Y)D_{Z}-\textbf{K}([X,Y],Z)\}=0$,
\end{description}
where $\mathfrak{S}_{X,Y,Z}$ denotes cyclic sum over the vector
fields $X, Y$ and $Z$.
\end{lem}

\par We terminate this section
by some concepts and results concerning the Klein-Grifone approach.
For more details, we refer to \cite{r21}, \cite{r22} and \cite{r27}.

 A semispray  on $M$ is a vector field $X$ on $TM$,
 $C^{\infty}$ on $\T M$, $C^{1}$ on $TM$, such that
$\rho\circ X = \overline{\eta}$. A semispray $X$ which is
homogeneous of degree $2$ in the directional argument
($[\mathcal{C},X]= X $) is called a spray.

\begin{prop}{\em{\cite{r27}}}\label{spray} Let $(M,L)$ be a Finsler manifold. The vector field
$G$ on $TM$ defined by $i_{G}\Omega =-dE$ is a spray, where
 $E:=\frac{1}{2}L^{2}$ is the energy function and $\Omega:=dd_{J}E$.
 Such a spray is called the canonical spray.
 \end{prop}

A nonlinear connection on $M$ is a vector $1$-form $\Gamma$ on $TM$,
$C^{\infty}$ on $\T M$, $C^{0}$ on $TM$, such that
$$J \Gamma=J, \quad\quad \Gamma J=-J .$$
The horizontal and vertical projectors $h_{\Gamma}$\,  and
$v_{\Gamma}$ associated with $\Gamma$ are defined by
   $h_{\Gamma}:=\frac{1}{2} (I+\Gamma),\, v_{\Gamma}:=\frac{1}{2}
 (I-\Gamma).$
Thus $\Gamma$ gives rise to the direct sum decomposition $T\T M=
H(\T M)\oplus V(\T M)$, where $H(\T M):=Im \, h_{\Gamma} = Ker\,
v_{\Gamma} $, $\,\,V(\T M):= Im \, v_{\Gamma}=Ker \, h_{\Gamma}$. We
have $ J\circ h_{\Gamma} =J, \,\,\, h_{\Gamma}\circ J=0, \,\,\,
J\circ v_{\Gamma}=0, \,\,\, v_{\Gamma}\circ J=J.$  A nonlinear
connection $\Gamma$ is homogeneous if $[\mathcal{C},\Gamma]=0$. The
torsion $t$ of a nonlinear connection $\Gamma$ is the vector
$2$-form  on $TM$ defined by $t:=\frac{1}{2} [J,\Gamma]$. The
curvature of $\Gamma$ is the vector $2$-form on $TM$ defined by
$\mathfrak{R}:=-\frac{1}{2}[h_{\Gamma},h_{\Gamma}]$. A nonlinear
connection $\Gamma$ is said to be conservative if
$d_{h_{\Gamma}}\,E=0$. With any given nonlinear connection $\Gamma$,
one can associate a semispray $S$ which is horizontal with respect
to $\Gamma$, namely, $S=h_{\Gamma}S'$, where $S'$ is an arbitrary
semispray. Moreover, if $\Gamma$ is homogeneous, then its associated
semispray is a spray.

\begin{thm} \label{th.9a} {\em{\cite{r22}}} On a Finsler manifold $(M,L)$, there exists a unique
conservative homogenous nonlinear  connection  with zero torsion. It
is given by\,{\em:} \vspace{-0.3cm} $$\Gamma =
[J,G],\vspace{-0.3cm}$$ where $G$ is the canonical spray.\\
 Such a nonlinear connection is called the canonical connection, or the Barthel connection, associated with $(M,L)$.
\end{thm}
\par
It should be noted that the semispray associated with the Barthel
connection is a spray, which is the canonical spray.


\Section{Cartan and Berwald connections in PB-formalism}

In this section, we review the main results obtained in \cite{r92}.
It concerns the spray and nonlinear connection associated with a
regular connection on $\pi^{-1}(TM)$  and the existence and
uniqueness theorems of Cartan and Berwald connections in Finsler
geometry. From these results, the relationships between the
curvature tensors associated with Berwald and Cartan connections are
investigated. \vspace{5pt}

\begin{defn} Let $D$ be a regular connection on
$\pi^{-1}(TM)$ with horizontal map $\beta$.\\
\textbf{--} The semispray  $S=\beta\circ\overline{\eta}$ will be called the semispray associated with $D$.\\
\textbf{--} The nonlinear connection $\Gamma=2\beta\circ\rho-I$ will
be called the nonlinear connection associated with $D$.
\end{defn}

\begin{prop}\label{eqv.}Let $(M,L)$ be a Finsler manifold. Let  ${D}$ be a regular connection  on
$\pi^{-1}(TM)$ whose connection map is $K$ and whose horizontal map
is $\beta$. Then, the following assertions are equivalent:
\vspace{-0.2cm}
 \begin{description}
    \item[(a)] The (h)hv-torsion  ${T}$ of ${D}$ has the property that
    ${T}( \overline{X},\overline{\eta})=0$,

    \item[(b)] $K=\gamma^{-1}$ on  $V(TM)$,

    \item[(c)] $\widetilde{\Gamma} :=\beta\circ\rho - \gamma\circ K$ is a nonlinear
    connection on $M$.\vspace{-0.2cm}
 \end{description}
 \par Consequently,  if any one of the above assertions
    holds, then $\widetilde{\Gamma}$
    coincides with the nonlinear connection associated with $D$: $\widetilde{\Gamma}=\Gamma=
    2\beta\circ\rho-I$, and in  this case $h_{\Gamma}=h_{D}=\beta\circ\rho$
    and  $\,v_{\Gamma}=v_{D}=\gamma\circ K$.
\end{prop}

\begin{lem}\label{bracket2} Let $D$ be a regular connection on $\p$
whose (h)hv-torsion tensor $T$ has the property that
$\,T(\overline{X}, \overline{\eta})=0$. Then, we
have{\em:}\vspace{-0.2cm}
   \begin{description}
 \item[(a)] $[\beta \overline{X},\beta \overline{Y}]=
     \gamma\widehat{R}(\overline{X},\overline{Y})
     + \beta(D_{\beta \overline{X}}\overline{Y}-
     D_{\beta \overline{Y}}\overline{X}-Q(\overline{X},\overline{Y})),$

    \item[(b)] $[\gamma \overline{X},\beta \overline{Y}]=-
     \gamma(\widehat{P}(\overline{Y},\overline{X})+D_
     {\beta \overline{Y}}\overline{X})
     +\beta( D_{\gamma \overline{X}}\overline{Y}-T(\overline{X},\overline{Y})),$

   \item[(c)] $[\gamma \overline{X},\gamma \overline{Y}]=
     \gamma(D_{\gamma \overline{X}}\overline{Y}-
     D_{\gamma \overline{Y}}\overline{X}+\widehat{S}(\overline{X},\overline{Y}))$.

     \end{description}
\end{lem}

 The following theorem guarantees the existence  and  uniqueness of the  Cartan connection.
 This is the Finsler analogue of the fundamental
theorem of Riemannian geometry.\vspace{-0.2cm}
\begin{thm}\label{le.cc3} Let $(M,L)$ be a Finsler
manifold and  $g$ the Finsler metric defined by $L$. There exists a
unique regular connection $\nabla$ on $\pi^{-1}(TM)$ such that
\begin{description}
  \item[(a)]  $\nabla$ is  metric\,{\em:} $\nabla g=0$,

  \item[(b)] The (h)h-torsion of $\nabla$ vanishes\,{\em:} $Q=0
  $,
  \item[(c)] The (h)hv-torsion $T$ of $\nabla$\, satisfies\,\emph{:}
   $g(T(\overline{X},\overline{Y}), \overline{Z})=g(T(\overline{X},\overline{Z}),\overline{Y})$.\vspace{-0.2cm}
\end{description}
\par
 Such a connection is called the Cartan
connection associated with  the Finsler manifold $(M,L)$.
\par
 This connection  is uniquely determined by  the  relations\,\emph{:}
   \begin{description}
     \item[(i)]  $ 2g(\nabla _{\gamma\overline{X}}
      \overline{Y}, \overline{Z})
      =\gamma\overline{X}\cdot g( \overline{Y},\overline{Z})+
     g( \overline{Y},\rho [\beta\overline{Z},\gamma\overline{X}])
     +g(\overline{Z},\rho [\gamma\overline{X},\beta \overline{Y}])$.

     \item[(ii)]  $ 2g(\nabla _{\beta\overline{X}}\rho Y,\rho Z)
     = \beta\overline{X}\cdot g( \overline{Y},
                \overline{Z})+
                \beta\overline{Y}\cdot g( \overline{Z},\overline{ X})
                -\beta\overline{Z}\cdot g( \overline{X},
                \overline{Y})$\\
  $ { \qquad \qquad \qquad \quad }
  -g( \overline{X},\rho [\beta\overline{Y},\beta\overline{Z}])+g( \overline{Y},\rho
[\beta\overline{Z},\beta\overline{X}])+g( \overline{Z},\rho
    [\beta\overline{X},\beta\overline{Y}])$.
 \end{description}
\end{thm}

Concerning the nonlinear connection associated with the Cartan
connection, we have \,:\vspace{-0.2cm}
\begin{thm} \label{c.ba.} Let $\nabla$ be the Cartan
connection. The nonlinear connection $\Gamma$ associated with
$\nabla$ coincides with the Barthel connection{\em:} $\Gamma=[J,G]$.
\end{thm}


Concerning  the existence and uniqueness of the Berwald connection,
we have\,:\vspace{-0.2cm}
\begin{thm}\label{bth2.h2} Let $(M,L)$ be a Finsler manifold. There exists a
unique regular connection ${{D}}^{\circ}$ on $\pi^{-1}(TM)$ such
that
\begin{description}
 \item[(a)] $D^{\circ}_{h^{\circ}X}L=0$,
  \item[(b)]   ${{D}}^{\circ}$ is torsion-free\,{\em:} ${\textbf{T}}^{\circ}=0 $,
  \item[(c)]The (v)hv-torsion tensor $\widehat{P^{\circ}}$ of ${D}^{\circ}$ vanishes\,\emph{:}
   $\widehat{P^{\circ}}(\overline{X},\overline{Y})= 0$.
  \end{description}
  \par Such a connection is called the Berwald
  connection associated with the Finsler manifold $(M,L)$.
\end{thm}

\begin{prop}\label{rem.1}
The semispray associated with the Berwald connection is a spray
which coincides with the canonical spray. Moreover, the nonlinear
connection associated with the Berwald connection coincides with the
Barthel connection.
\end{prop}

\begin{thm} \label{th.5} The Berwald connection $D^{\circ}$ is explicitly expressed in
terms of the Cartan connection $\nabla$ in the form\,\textsc{:}
 \vspace{-0.1cm}
  \begin{equation}\label{b10}
     {{D}}^{\circ}_{X}\overline{Y} = \nabla _{X}\overline{Y}
+{\widehat{P}}(\rho X,\overline{Y}) -T(K X,\overline{Y}).
\vspace{-0.2cm}
\end{equation}
In particular, we have\vspace{-0.1cm}
\begin{description}
  \item[(a)] $ {{D}}^{\circ}_{\gamma \overline{X}}\overline{Y}=\nabla _{\gamma
  \overline{X}}\overline{Y}-T(\overline{X},\overline{Y})$.

 \item[(b)] $ {{D}}^{\circ}_{\beta \overline{X}}\overline{Y}=\nabla _{\beta
  \overline{X}}\overline{Y}+\widehat{P}(\overline{X},\overline{Y}).$
\end{description}
\end{thm}

In view of the above theorem, we have\,: \vspace{-0.11cm}
\begin{prop} The Berwald connection $D^{\circ}$ has the properties:
\begin{description}
\item[(a)] $({{D}}^{\circ}_{\gamma
\overline{X}}\,g)(\overline{Y},\overline{Z}) =2 g(T(\overline{X},
\overline{Y}),\overline{Z})$.

\item[(b)] $({{D}}^{\circ}_{\beta \overline{X}}\,g)(\overline{Y},\overline{Z}) =-2
  g(\widehat{P}(\overline{X},\overline{Y}),\overline{Z})$.

\item[(c)] $ {{D}}^{\circ}_{G}\,g=0$.
\end{description}
Consequently,\\
\textbf{--} a Finsler manifold $(M,L)$ is Riemannian if and only if
$\,{{D}}^{\circ}_{\gamma \overline{X}}\,g=0$.\\
\textbf{--} a Finsler manifold $(M,L)$ is Landsbergian if and only
if $\,{{D}}^{\circ}_{\beta \overline{X}}\,g=0$.
\end{prop}

\begin{prop}\label{b7.pp.7} The Curvature tensor $\textbf{K}$ of the Cartan connection
$\nabla$ and The Curvature tensor $\textbf{K}^{\circ}$  of the
 Berwald connection  ${D}^{\circ}$ are
related by\,\textsc{:}\vspace{-0.1cm}
\begin{eqnarray*}
  \textbf{{K}}^{\circ}(X,Y)\overline{Z} &=&\textbf{K}(X,Y)\overline{Z}-
  \widehat{P}(\textbf{T}(X,Y),\overline{Z}) -T(K[X,Y], \overline{Z})
  -\mathfrak{U}_{X,Y}\{(\nabla_{X}\widehat{P})(\rho Y,\overline{Z})
   \\
  &&-\nabla_{X}T(K Y, \overline{Z})+T(K Y,\nabla_{X}\overline{Z})
  +\widehat{P}(\rho X,\widehat{P}(\rho Y,\overline{Z}))-\widehat{P}(\rho X, T(K Y,\overline{Z}))\\
   &&+T(K X,T(K Y,\overline{Z}))-T(K X,\widehat{P}(\rho Y,\overline{Z}))\},\vspace{-0.2cm}
\end{eqnarray*}
where \,$\mathfrak{U}_{X,Y}A(X,Y)=A(X,Y)-A(Y,X)$.\vspace{5pt}\\
 In particular, we
have\vspace{-0.1cm}
\begin{description}

  \item[(a)]${{S}}^{\circ}(\overline{X},\overline{Y})\overline{Z}= 0$.

  \item[(b)]${{P}}^{\circ}(\overline{X},\overline{Y})\overline{Z}= P(\overline{X},\overline{Y})\overline{Z}+
  (\nabla_{\gamma \overline{Y}}\widehat{P})(\overline{X},\overline{Z}) +
   \widehat{P}(T(\overline{Y},\overline{X}),\overline{Z})
   +\widehat{P}(\overline{X},T(\overline{Y},\overline{Z}))+$\\
${\qquad\qquad\ \ \ }
+(\nabla_{\beta\overline{X}}T)(\overline{Y},\overline{Z})
-T(\overline{Y}, \widehat{P}(\overline{X},\overline{Z}))-T(
\widehat{P}(\overline{X}, \overline{Y}), \overline{Z})$.
  \item[(c)]${{R}}^{\circ}(\overline{X},\overline{Y})\overline{Z}= R(\overline{X},\overline{Y})\overline{Z}- T(\widehat{R}(\overline{X},\overline{Y}),\overline{Z})-
  \mathfrak{U}_{\overline{X},\overline{Y}}\{(\nabla_{\beta \overline{X}}\widehat{P})(\overline{Y}, \overline{Z})
    +\widehat{P}(\overline{X},\widehat{P}(\overline{Y},\overline{Z}))\}.$

\end{description}
\end{prop}

\begin{cor}~\par\vspace{-0.2cm}
\begin{description}
    \item[(a)] $\widehat{S^{\circ}}(\overline{X}, \overline{Y})=0$,

    \item[(b)] $\widehat{P^{\circ}}(\overline{X}, \overline{Y} )=0$,
    \item[(c)]  $ \widehat{R^{\circ}}(\overline{X},\overline{Y})= \widehat{R}(\overline{X},\overline{Y})$.
\end{description}
\end{cor}

\Section{Chern (Rund) connection}

 In this section, we establish an  intrinsic
(coordinate-free) proof of the existence and  uniqueness theorem  of
 Chern (Rund) connection. Moreover, the relationships between
this connection and the Cartan  and Berwald connections are
obtained.

 \vspace{7pt}
\par
We start with the following fundamental result. \vspace{-0.2cm}
\begin{thm} \label{th.a} Let $D^{\diamond}$ be a regular connection on $\pi^{-1}(TM)$
 with  connection map $K^{\diamond}$ such that
\begin{description}
  \item[(a)]  $(D^{\diamond}_{X}\,g)(\rho Y,\rho Z)=2g(T(K^{\diamond} X,\rho Y), \rho
  Z)$,

  \item[(b)]   $D^{\diamond}$ is torsion free\,\emph{:} ${\textbf{T}}\,^{\diamond}=0$.
\end{description}
Then, the nonlinear connection $\Gamma^{\diamond}$ associated with
$D^{\diamond}$ coincides with the Barthel connection{\,\em:}
$\Gamma=[J,G]$.
\end{thm}

To prove this theorem, we need the following three lemmas.968
\vspace{-0.2cm}
\begin{lem}\label{lem.a1} The hv-curvature tensor ${P}^{\diamond}$ of the  connection
$D^{\diamond}$ is symmetric with respect to  the first and  third
arguments\,{\em:}
    $$P^{\diamond}(\overline{X},\overline{Y})\overline{Z}=
    P^{\diamond}(\overline{Z},\overline{Y})\overline{X} \ \ \ for \ all
    \  \overline{X}, \overline{Y}, \overline{Z }\in \cp .$$
\end{lem}

\begin{proof}
    The proof follows  from Lemma \ref{cyclic}(a) by setting
    $X=\beta^{\diamond}\, \overline{X}$, $Y=\gamma \overline{Y}$ and $Z=\beta^{\diamond}\,\overline{Z}$, noting that
    $\rho\circ \gamma=0$, $\rho\circ\beta^{\diamond}\,=id_{\cp}$ and that $\textbf{T}^{\diamond}=0$.
\end{proof}

\begin{lem}\label{lem.a2}The hv-curvature tensor $P^{\diamond}$ of the  connection
$D^{\diamond}$ has the property that
$$P^{\diamond}(\overline{X},\overline{Y},\overline{Z},\overline{W})
    +P^{\diamond}(\overline{X},\overline{Y},\overline{W},\overline{Z})
    =2(D^{\diamond}_{\beta^{\diamond}\,{\overline{X}}}T)(\overline{Y},\overline{Z},\overline{W})
    -2T(\widehat{P^{\diamond}}(\overline{X},\overline{Y}),
    \overline{Z},\overline{W}).$$
\end{lem}

\begin{proof} We have
$$ X\cdot g(\overline{W},\overline{Z} )=(D^{\diamond}_{X}g)(\overline{W},\overline{Z})
+g(D^{\diamond}_{X}\overline{W},\overline{Z})+g(\overline{W},D^{\diamond}_{X}\overline{Z}).$$
From which, we obtain
\begin{eqnarray*}
   X\cdot ( Y\cdot g(\overline{W},\overline{Z} ))&=&
     X\cdot((D^{\diamond}_{Y}g)(\overline{W},\overline{Z}))
+ X\cdot g(D^{\diamond}_{Y}\overline{W},\overline{Z})+
 X\cdot g(\overline{W},D^{\diamond}_{Y}\overline{Z})\\
   &=&  X\cdot((D^{\diamond}_{Y}g)(\overline{W},\overline{Z}))
   + (D^{\diamond}_{X}g)(D^{\diamond}_{Y}\overline{W},\overline{Z})+
   (D^{\diamond}_{X}g)(\overline{W},D^{\diamond}_{Y}\overline{Z})
   \\
   &&+g(D^{\diamond}_{X}D^{\diamond}_{Y}\overline{W},\overline{Z})
   +g(D^{\diamond}_{Y}\overline{W},D^{\diamond}_{X}\overline{Z})
   +g(D^{\diamond}_{X}\overline{W},D^{\diamond}_{Y}\overline{Z})+\\
   && +g(\overline{W},D^{\diamond}_{X}D^{\diamond}_{Y}\overline{Z}),
  \end{eqnarray*}
with similar expression for $Y\cdot ( X\cdot
g(\overline{W},\overline{Z} ))$. Then,
\begin{eqnarray*}
   [X,Y]\cdot g(\overline{W},\overline{Z} )&=&(D^{\diamond}_{[X,Y]}g)(\overline{W},\overline{Z})
+g(D^{\diamond}_{[X,Y]}\overline{W},\overline{Z})+g(\overline{W},D^{\diamond}_{[X,Y]}\overline{Z})\\
&=& \mathfrak{U}_{X,Y}\{
X\cdot((D^{\diamond}_{Y}g)(\overline{W},\overline{Z}))
   + (D^{\diamond}_{X}g)(D^{\diamond}_{Y}\overline{W},\overline{Z})+
   (D^{\diamond}_{X}g)(\overline{W},D^{\diamond}_{Y}\overline{Z})\}
   \\
   &&+g([D^{\diamond}_{X},D^{\diamond}_{Y}]\overline{W},\overline{Z})
      +g(\overline{W},[D^{\diamond}_{X},D^{\diamond}_{Y}]\overline{Z}).
  \end{eqnarray*}
Hence,
\begin{equation}\label{000}
\left.
    \begin{array}{rcl}
  &&
  g(\textbf{K}(X,Y)\overline{Z},\overline{W})+g(\textbf{K}(X,Y)\overline{W},\overline{Z})=\\
    &=& 2\mathfrak{U}_{X,Y}\{
X\cdot g(T(K^{\diamond} Y,\overline{W}),\overline{Z})
   + g(T(K^{\diamond} X,D^{\diamond}_{Y}\overline{W}),\overline{Z}) \\
   &&+
   g(T(K^{\diamond} X,\overline{W}),D^{\diamond}_{Y}\overline{Z})\}
  -2g(T(K^{\diamond}[X,Y],\overline{W}),\overline{Z}).
 \end{array}
  \right.
\end{equation}
From which, by setting  $X=\beta^{\diamond}\, \overline{X}$ and
$Y=\gamma \overline{Y}$ into (\ref{000}) and using Lemma
\ref{bracket2}, we get
\begin{eqnarray*}
   P^{\diamond}(\overline{X},\overline{Y},\overline{Z},\overline{W})
   +P^{\diamond}(\overline{X},\overline{Y},\overline{W},\overline{Z})
&=&2\beta^{\diamond}\,\overline{X}\cdot T(
\overline{Y},\overline{W},\overline{Z})
   - 2T(\overline{Y},D^{\diamond}_{\beta^{\diamond}\,\overline{X}}\overline{W},\overline{Z}) \\
   &&- 2T(\overline{Y},\overline{W},D^{\diamond}_{\beta^{\diamond}\,\overline{X}}\overline{Z})
  -2T(\widehat{P^{\diamond}}(\overline{X},\overline{Y}),\overline{W},\overline{Z})-\\
  &&-2T( D^{\diamond}_{\beta^{\diamond}\,\overline{X}}\overline{Y},\overline{W},\overline{Z}).
  \end{eqnarray*}
Hence, the result follows.
\end{proof}

\begin{lem}\label{lem.a3}The (v)hv-torsion tensor $\widehat{P^{\diamond}}$  given
by\,:
 $$\widehat{P^{\diamond}}(\overline{X}, \overline{Y})=
   (D^{\diamond}_{\beta^{\diamond}\,{\overline{\eta}}}T)(\overline{X},\overline{Y}).$$
 Consequently,  $\widehat{P^{\diamond}}$ is symmetric and
$\widehat{P^{\diamond}}(\overline{X},\overline{\eta})=0$.
\end{lem}

\begin{proof}Firstly, one can easily show that\vspace{-0.2cm}
 \begin{equation}\label{eq}
   (D^{\diamond}_{\beta^{\diamond}\,{\overline{X}}}T)(\overline{Y},\overline{Z},\overline{W})=
g((D^{\diamond}_{\beta^{\diamond}\,{\overline{X}}}T)(\overline{Y},\overline{Z}),\overline{W}).\vspace{-0.2cm}
 \end{equation}
Cyclic permutation on $\overline{X}, \overline{Z}, \overline{W}$ in
the formula of Lemma \ref{lem.a2} yields three equations. Adding two
of these equations and subtracting the third gives\vspace{-0.2cm}
 \begin{equation}\label{01}
\left.
    \begin{array}{rcl}
   P^{\diamond}(\overline{X}, \overline{Y}, \overline{Z}, \overline{W}) &=&
   (D_{\beta^{\diamond}\, \overline{X}}T)(\overline{Y}, \overline{Z}, \overline{W})
   +(D_{\beta^{\diamond}\, \overline{Z}}T)(\overline{Y}, \overline{W}, \overline{X})
   -(D_{\beta^{\diamond}\,\overline{W}}T)(\overline{Y}, \overline{X}, \overline{Z}) \\
    &&+T(\widehat{P^{\diamond}}(\overline{W},\overline{Y}), \overline{X},\overline{Z})
    -T(\widehat{P^{\diamond}}(\overline{X},\overline{Y}), \overline{Z},\overline{W})
    -T(\widehat{P^{\diamond}}(\overline{Z},\overline{Y}),
    \overline{W},\overline{X}).
\end{array}
  \right.
\end{equation}

Setting $\overline{X}=\overline{Z}=\overline{\eta}$ in (\ref{01}),
taking into account the properties of the
    (h)hv-torsion $T$  and the fact that  $K^{\diamond}\circ \beta^{\diamond}\,=0$,
     we obtain\vspace{-0.2cm}
     \begin{equation}\label{0}
       \widehat{P^{\diamond}}(\eta, \overline{Y})=0 ,\ \forall \  \overline{Y} \in \cp.\vspace{-0.2cm}
     \end{equation}
Again,  setting $ \overline{Z}=\overline{\eta}$ in (\ref{01}) and
now using (\ref{0}) and (\ref{eq}), we conclude that\vspace{-0.2cm}
\begin{equation}\label{02}
   \widehat{P^{\diamond}}(\overline{X}, \overline{Y})=
   (D^{\diamond}_{\beta^{\diamond}\,{\overline{\eta}}}T)(\overline{X},\overline{Y}).\vspace{-0.2cm}
\end{equation}
The symmetry of $ \widehat{P^{\diamond}}$ follows then from
(\ref{02}) and  the symmetry of $T$.
\end{proof}

\noindent\textit{\textbf{Proof of Theorem \ref{th.a}}}\,:
\par
 As the (h)hv-torsion tensor $T^{\diamond}$  of the connection
$D^{\diamond}$ vanishes, then
$T^{\diamond}(\overline{X},\overline{\eta})=0$. Consequently,  by
Proposition \ref{eqv.}, it follows that the associated nonlinear
connection  $\Gamma^{\diamond} $ has the form \vspace{-0.2cm}
$$\Gamma^{\diamond} :=\beta^{\diamond}\, o \rho - \gamma o K^{\diamond} \vspace{-0.2cm}.$$
\par
 Now, we prove that the connection $\Gamma^{\diamond}$ satisfies the
 following properties:
 \par
 \vspace{5pt}
 \noindent\textit{\textbf{$\Gamma^{\diamond} $ is conservative:}}
   \noindent  ${d}_{{h^{\diamond}}X}E=0,$ for all $X\in\cpp$:\\
    In fact,  by condition (a), $D^{\diamond}_{h^{\diamond}X} g=0$. Taking this,
    together with the identities  $2E=g(\overline{\eta},
  \overline{\eta})$ and $ K^{\diamond}\circ h^{\diamond}=0$ into account,
    we get
     \begin{eqnarray*}
      {d}_{{h^{\diamond}}X}E&=& h^{\diamond}X\cdot E=
      \frac{1}{2} h^{\diamond}X \cdot g(\overline{\eta}, \overline{\eta})
      = g( D^{\diamond}_{h^{\diamond}X}\overline{\eta},\overline{\eta})
      =g( K^{\diamond}({h^{\diamond}X}),\overline{\eta})=0.
     \end{eqnarray*}
 \vspace{5pt}
   \noindent \textit{\textbf{$\Gamma^{\diamond}$ is
homogenous $(\,[\mathcal{C},\Gamma^{\diamond}]=0)$:}}
\par
It is easy to show that
  $$
    [\mathcal{C},v^{\diamond}]X=-v^{\diamond}[\mathcal{C},h^{\diamond}X]
  $$
As $v^{\diamond}=\gamma\circ K^{\diamond}$,
$h^{\diamond}=\beta^{\diamond}\,\circ\rho$ and
$\gamma\circ\overline{\eta}=\mathcal{C}$, then

$$[\mathcal{C},v^{\diamond}]X=-(\gamma \circ K^{\diamond})[\gamma \overline{\eta},\beta^{\diamond}\,\rho X]$$
Using Lemma \ref{bracket2}, we obtain
\begin{eqnarray*}
    [\mathcal{C},v^{\diamond}]X
    &=&-(\gamma\circ K^{\diamond})\{-
     \gamma(\widehat{P^{\diamond}}(\rho X,\overline{\eta})+D^{\diamond}_
     {\beta^{\diamond}\, \rho X}\overline{\eta})
     +\beta^{\diamond}\, D^{\diamond}_{\gamma \overline{\eta}}\rho X\}\\
     &=&
     \gamma\widehat{P^{\diamond}}(\rho X,\overline{\eta}), \,\,
     as \,\,  \ K^{\diamond}\circ
     \beta^{\diamond}=0\, \ and \ K^{\diamond}\circ\gamma=id_{\cp}.
  \end{eqnarray*}
Finally, by Lemma \ref{lem.a3}, $ \widehat{P^{\diamond}}(
\overline{X},\overline{\eta})=0 $ and consequently, $
[\mathcal{C},\Gamma^{\diamond}]= -2[\mathcal{C},v^{\diamond}]=0$.
\par
 \vspace{5pt}
 \noindent
\textit{\textbf{$\Gamma^{\diamond}$ is torsion-free
$(\,[J,\Gamma^{\diamond}]=0)$:}}
\begin{eqnarray*}
  [J,v^{\diamond}](X,Y) &=& [JX,v^{\diamond}Y]+ [v^{\diamond}X,JY]+v^{\diamond}J[X,Y]+Jv^{\diamond}[X,Y]\\
   && -J[v^{\diamond}X,Y]-J[X,v^{\diamond}Y]-v^{\diamond}[JX,Y]-v^{\diamond}[X,JY].
\end{eqnarray*}
As $J\circ v^{\diamond}=0$, $v^{\diamond}\circ J=J$ and the vertical
distribution is completely integrable, we get
\begin{eqnarray*}
[J,v^{\diamond}](X,Y) &=&J[h^{\diamond}X,h^{\diamond}Y]-v^{\diamond}[JX,h^{\diamond}Y]-v^{\diamond}[h^{\diamond}X,JY]\\
   &=&J[\beta^{\diamond}\, \rho X,\beta^{\diamond}\, \rho Y]-
   v^{\diamond}[\gamma \rho X,\beta^{\diamond}\, \rho Y]+v^{\diamond}[\gamma \rho Y,\beta^{\diamond}\, \rho
   X].
\end{eqnarray*}
From which, together with Lemma \ref{bracket2}, we obtain
  \begin{eqnarray*}
  [J,v^{\diamond}](X,Y)&=&J\{\gamma\widehat{R^{\diamond}}(\rho{X},\rho{Y})
     + \beta^{\diamond}\,(D^{\diamond}_{h^{\diamond}{X}}\rho{Y}-
     D^{\diamond}_{h^{\diamond}{Y}}\rho{X})\}\\
    &&-(\gamma\circ K^{\diamond})\{-
     \gamma(\widehat{P^{\diamond}}(\rho{Y},\rho{X})+D^{\diamond}_
     {h{Y}}\rho{X})
     +\beta^{\diamond}\,( D^{\diamond}_{J{X}}\rho{Y})\}\\
     &&+(\gamma\circ K^{\diamond})\{-
     \gamma(\widehat{P^{\diamond}}(\rho{X},\rho{Y})+D^{\diamond}_
     {h^{\diamond}{X}}\rho{Y})
     +\beta^{\diamond}\,( D^{\diamond}_{J{Y}}\rho{X})\}
     \end{eqnarray*}
Noting that $\,J\circ\gamma=0$, $\,J\circ\beta^{\diamond}\,=\gamma$,
$\,K^{\diamond}\circ\gamma=id_{\cp}$,
$\,K^{\diamond}\,\circ\beta^{\diamond}\,=0$
        and that
    $ \widehat{P^{\diamond}}$ is symmetric (by Lemma \ref{lem.a3}), it follows that $[J,v^{\diamond}]=0$.
    From which
    $t:=\frac{1}{2}[J,\Gamma^{\diamond}]=-[J,v^{\diamond}]=0$.
   \par From the above consideration, $\Gamma^{\diamond}=\beta^{\diamond}\,\circ\rho -\gamma\circ K^{\diamond}$
   is a conservative torsion-free homogenous
   nonlinear connection. By the uniqueness of the Barthel connection
    (Theorem \ref{th.9a}), it follows that  $\Gamma^{\diamond}$ coincides with the Barthel
   connection $[J,G]$.\ \  \ \ \ \ $\Box$

\vspace{7pt}

 In view of  Theorem \ref{th.a}, Theorem \ref{c.ba.} and Proposition \ref{rem.1}, we have\,:
 \vspace{-0.2cm}
\begin{cor}\label{rcor2} The nonlinear connection  associated with the  connection $D^{\diamond}$
is the same as the nonlinear connection associated with  the Cartan
connection $\nabla$ or  the Berwald connection $D^{\circ}$, which  coincides
 with  the Barthel connection $[J,G]$.\\
Consequently, \  $h^{\diamond}=h^{\circ}=h=\beta \circ \rho$,
$v^{\diamond}=v^{\circ}=v=\gamma \circ K$ and hence
  $\beta^{\diamond}\,=\beta^{\circ}=\beta$,
$K^{\diamond}=K^{\circ}=K$.
\end{cor}

Now, we are in a position to announce the main result of this
section.  \vspace{-0.2cm}
\begin{thm} \label{th.r1} Let $(M,L)$ be a Finsler manifold and $g$ the Finsler metric
defined by $L$. There exists a unique regular connection
$D^{\diamond}$ on $\pi^{-1}(TM)$ such that
\begin{description}
  \item[(a)]  $(D^{\diamond}_{X}\,g)(\rho Y,\rho Z)=2g(T(K^{\diamond} X,\rho Y), \rho
  Z)$,

  \item[(b)]   $D^{\diamond}$ is torsion free\,\emph{:} $\textbf{T}^{\diamond}=0$,
\end{description}
\par
\noindent where $T$ is the (h)hv-torsion of the Cartan connection
and
$K^{\diamond}$ is the connection map of $D^{\diamond}$.\\
This connection is called the Chern connection associated with
$(M,L)$.
\end{thm}

\begin{proof}
  First we prove the {\it \textbf{uniqueness}}. If $D^{\diamond}$ is
a non-metric linear connection on $\pi^{-1}(T M)$ with
 nonzero torsion $\textbf{T}^{\diamond}$, then
 $D^{\diamond}$  is completely  determined by\,:
\vspace{-0.2cm}
\begin{equation}\label{eq.r}
     \left.
    \begin{array}{rcl}
    2g(D^{\diamond} _{X}\rho Y,\rho Z)& =& X\cdot g(\rho Y,\rho
                Z)+ Y\cdot g(\rho Z,\rho X)-Z\cdot g(\rho X,\rho
                Y) \\
        & &-g(\rho X,\textbf{T}^{\diamond}(Y,Z))+g(\rho Y,\textbf{T}^{\diamond}(Z,X))
        +g(\rho Z,\textbf{T}^{\diamond}(X,Y)) \\
        & &-g(\rho X,\rho [Y,Z])+g(\rho Y,\rho [Z,X])+g(\rho Z,\rho [X,Y])\\
        & &- (D^{\diamond} _{X}\, g)(\rho Y,\rho
                Z)-(D^{\diamond} _{Y}\, g)(\rho Z,\rho X)+(D^{\diamond} _{Z}\,g)(\rho X,\rho
                Y).\vspace{-0.2cm}
   \end{array}
  \right\}
 \end{equation}
 for all $X,Y,Z\in\cpp$.  The connection $D^{\diamond}$ being
regular, let $h^{\diamond}$ and $v^{\diamond}$ be its  horizontal
and vertical projectors:
 $h^{\diamond}=\beta^{\diamond}\,\circ\rho,\, v^{\diamond}=I-\beta^{\diamond}\,\circ\rho $.
  From  Corollary \ref{rcor2}, we have
  $h^{\diamond}=h=\beta \circ \rho$ and
$v^{\diamond}=v=\gamma \circ K$.
 \par
 Now, by replacing $X, Y, Z$ by ${h}X, {h}Y, {h}Z $
in (\ref{eq.r})
 and using conditions (a) and (b), taking into account
  the fact that $\rho \circ {h}=\rho$, $K\circ  h=0$, we get
\begin{equation*}
    \left.
    \begin{array}{rcl}
     2g(D^{\diamond} _{hX}\rho Y,\rho Z) & = & hX\cdot g(\rho Y,\rho
                Z)+ hY\cdot g(\rho Z,\rho X)-hZ\cdot g(\rho X,\rho
                Y) \\
        & & -g(\rho X,\rho [hY,hZ])+g(\rho Y,\rho
[hZ,hX])+g(\rho Z,\rho
    [hX,hY]).
    \end{array}
  \right.
\end{equation*}
From which, together with Theorem \ref{le.cc3}(ii), we get
\begin{equation}\label{eq.rt5}
D^{\diamond} _{hX}\rho Y=\nabla _{hX}\rho Y.
\end{equation}

Similarly, by replacing $X, Y, Z$ by $vX, vY, vZ $ in (\ref{eq.r}),
where $vX=\gamma \overline{X}$ for some $\overline{X}\in \cp$  and
using conditions (a) and (b), noting that $\rho \circ v=0$, we have
     \begin{equation*}\label{eq.rt6}
    \left.
    \begin{array}{rcl}
              2g(D^{\diamond} _{vX}\rho Y,\rho Z) &=&vX\cdot g(\rho Y,\rho Z)+
     g(\rho Y,\rho [{h}Z,{v}X])+g(\rho Z,\rho [{v}X,hY]).\\
          &&- (D^{\diamond} _{vX}\, g)(\rho Y,\rho
                Z)\\
                &=&\gamma \overline{X}\cdot g(\rho Y,\rho Z)+
     g(\rho Y,\rho [{h}Z,\gamma \overline{X}])+g(\rho Z,\rho [\gamma \overline{X},{h}Y]).\\
          &&- 2g(T( \overline{X}, \rho Y), \rho Z).
    \end{array}
  \right.
\end{equation*}
From which, together with Theorem \ref{le.cc3}(i), we have
\begin{equation}\label{eq.rt6}
    D^{\diamond} _{{v}X}\rho Y=\nabla _{{v}X}\rho Y-T(KX,\rho Y)=\rho [{v}X,{h}Y] .
\end{equation}
Now, from  Equations (\ref{eq.rt5}) and (\ref{eq.rt6}), we get
\begin{equation}\label{03}
D^{\diamond}_{X}\overline{Y} = \nabla _{X}\overline{Y} -
T(KX,\overline{Y}).
\end{equation}
 Hence $D^{\diamond}_{X}\rho Y$ is
uniquely determined by  Equations (\ref{03}).
\par
\vspace{4pt}
 To show the {\it \textbf{existence}}, we define $D^{\diamond}$ by
the requirement that (\ref{03}) holds for all $X, Y \in \x$. Now, we
need to prove that the connection $D^{\diamond}$
 satisfies the following properties\,:
 \par
 \vspace{3pt}
\noindent\textit{\textbf{ $D^{\diamond}$ satisfies condition
(a)}}\,: By using (\ref{03}) and the properties of the Cartan
connection $\nabla$, we get\vspace{-0.2cm}
 \begin{eqnarray*}
   (D^{\diamond}_{{X}}\,g)(\overline{Y},\overline{Z}) &=&
   D^{\diamond}_{{X}}\,g(\overline{Y},\overline{Z})-
g(D^{\diamond}_{{X}}\overline{Y},\overline{Z})-
g(\overline{Y},D^{\diamond}_{{X}}\overline{Z})\\
&=&
   X \cdot g(\overline{Y},\overline{Z})-
g(\nabla _{X}\overline{Y}-T(K X, \overline{Y}),\overline{Z})-
 g(\overline{Y},\nabla _{{X}}\overline{Z}-T(K X,\overline{Z}))\\
   &=&(\nabla _{{X}}g)(\overline{Y},\overline{Z})+
  g(T(K {X},\overline{Y}),\overline{Z})+
  g(\overline{Y},T(K {X},\overline{Z}))\\
  &=&2 g(T(K {X},\overline{Y}),\overline{Z})=2 g(T(K^{\diamond}{X},\overline{Y}),\overline{Z}).\vspace{-0.1cm}
\end{eqnarray*}

\par
 \vspace{3pt}
\noindent\textit{\textbf{ $D^{\diamond}$ satisfies condition
\emph{(}b\emph{)}}}\,:
 Again using (\ref{03})and the properties of the $\nabla$, we get  \vspace{-0.2cm}
\begin{eqnarray*}
  \textbf{T}^{\diamond}(X,Y)&=&D^{\diamond}_{{X}}\rho Y -D^{\diamond}_{Y}\rho
  X-\rho[X,Y]\\
  &=& \nabla _{X}\rho{Y} -
T(KX,\rho{Y})-\nabla _{Y}\rho{X} + T(KY,\rho{X})-\rho[X,Y]\\
&=&\textbf{T}(X,Y)-\textbf{T}(vX,hY)+\textbf{T}(vY,hX)=0.\vspace{-0.2cm}
\end{eqnarray*}
 This completes the proof.
\end{proof}

It is to be  noted  that  when we localize the above result, the {
\it{local expressions}} obtained coincide with the classical
expressions found in \cite{r91}, \cite{r23}, \cite{r93}...etc. (c.f.
Appendix 2).
\par
In view of the above theorem, taking  Theorem \ref{th.5} into
account, we have\vspace{-0.2cm}
\begin{thm}\label{th.r5}The Chern connection $D^{\diamond}$
is given in terms of the Cartan connection $\nabla$ \emph{(}or the
Berwald connection $D^{\circ}$\emph{)} by\,:
   $$ D^{\diamond}_{X}\overline{Y} = \nabla _{X}\overline{Y}
- T(KX,\overline{Y})= D^{\circ} _{X}\overline{Y} -{\widehat{P}}(\rho
X, \overline{Y}). \vspace{-0.2cm}$$ In particular, we have
\begin{description}
  \item[(a)] $ D^{\diamond}_{\gamma \overline{X}}\overline{Y}=\nabla _{\gamma
  \overline{X}}\overline{Y}-T(\overline{X},\overline{Y})=D^{\circ} _{\gamma
  \overline{X}}\overline{Y}$.

 \item[(b)] $ D^{\diamond}_{\beta \overline{X}}\overline{Y}=\nabla _{\beta
  \overline{X}}\overline{Y}=D^{\circ} _{\beta
  \overline{X}}\overline{Y}-{\widehat{P}}(\overline{X}, \overline{Y}).$
\end{description}
\end{thm}

\begin{cor}\label{rcor1}~\par\vspace{-0.2cm}
\begin{description}
\item[(a)] $ (D^{\diamond}_{\beta^{\diamond}\,\overline{X}}g)(\overline{Y},\overline{Z}) =
    0$.
    \item[(b)]  $(D^{\diamond}_{\gamma \overline{X}}g)(\overline{Y},\overline{Z})
=2 g(T(\overline{X}, \overline{Y}),\overline{Z})$.\vspace{-0.2cm}
\end{description}
Consequently, Finsler manifold $(M,L)$ is Riemannian if, and only
if,
  $D^{\diamond}_{\gamma \overline{X}}g=0$.
\end{cor}

Concerning the curvature tensor of the Chern connection
$D^{\diamond}$, we have\,:\vspace{-0.2cm}
\begin{prop} The curvature tensor of the Chern connection  $D^{\diamond}$
is given, in terms of the  curvature tensor of the Cartan connection
$\nabla$, by\,:\vspace{-0.1cm}
\begin{equation*}
\left.
    \begin{array}{rcl}
   \textbf{K}^{\,\diamond}(X,Y)\overline{Z} &=&\textbf{K}(X,Y)\overline{Z}-T(K[X,Y],
  \overline{Z})+\mathfrak{U}_{X,Y}\{\nabla_{X}T(K Y,\overline{Z})\\
  && -T(K Y,\nabla_{X}\overline{Z}) -T(K X,T(K Y,\overline{Z}))\}\vspace{-0.1cm}.
\end{array}
  \right.
\end{equation*}
In particular, we have\vspace{-0.1cm}
\begin{description}

  \item[(a)]$S^{\diamond}(\overline{X},\overline{Y})\overline{Z}= 0$.

  \item[(b)]$P^{\diamond}(\overline{X},\overline{Y})\overline{Z}=
  P(\overline{X},\overline{Y})\overline{Z}-T(\widehat{P}(\overline{X},\overline{Y}), \overline{Z})
  +(\nabla_{\beta \overline{X}}T)(\overline{Y},\overline{Z})$.

  \item[(c)]$R^{\diamond}(\overline{X},\overline{Y})\overline{Z}=
  R(\overline{X},\overline{Y})\overline{Z}-T(\widehat{R}(\overline{X},\overline{Y}),
  \overline{Z})$.
\end{description}
\end{prop}

\begin{cor}~\par\vspace{-0.2cm}
\begin{description}
    \item[(a)] $\widehat{S^{\diamond}}(\overline{X}, \overline{Y})=0$,

    \item[(b)] $\widehat{P^{\diamond}}(\overline{X}, \overline{Y} )=\widehat{P}(\overline{X}, \overline{Y}
    )$.
    \item[(c)]  $ \widehat{R^{\diamond}}(\overline{X},\overline{Y})= \widehat{R}(\overline{X},\overline{Y})$.
\end{description}
\end{cor}


\Section{Hashiguchi connection}

 In this section we establish an intrinsic proof of the
existence and  uniqueness theorem of the  Hashiguchi connection.
Moreover, the relationship between this connection and the Cartan
connection $\nabla$ is obtained.

\begin{thm} \label{0th.h2} Let $(M,L)$ be a Finsler manifold and  $g$  the Finsler metric
defined by $L$. Let $D^{*}$ be a regular connection on
$\pi^{-1}(TM)$ such that
\begin{description}

  \item[(a)]  ${D}^{*}$ is vertically  metric\,{\em:} ${D}^{*}_{\gamma \overline{X}}\, g=0$,

  \item[(b)]  The (h)hv-torsion tensor $T^{*}$ of ${D}^{*}$ satisfies\,{\em:}
  $g({T}^{*}(\overline{X},\overline{Y}), \overline{Z})=
  g({T}^{*}(\overline{X},\overline{Z}),\overline{Y})$,

  \item[(c)]   The (h)h-torsion of ${D}^{*}$ vanishes\,{\em:} $Q^{*}=0
  $,

  \item[(d)]The (v)hv-torsion  of ${D}^{*}$
  vanishes\,{\em:} $\widehat{{P}^{*}}=0$
 \item[(e)] $D^{*}_{h^{*}X}L=0$, $h^{*}$ being the horizontal
    projector of $D^{*}$.
  \end{description}
 Then, the nonlinear
connection $\Gamma^{*}$ associated with the   connection ${D}^{*}$
coincides with the Barthel connection {\em:} $\Gamma=[J,G]$.
\end{thm}

To prove this theorem, we need the following lemma\,:
\vspace{-0.2cm}
\begin{lem} \label{lem.ah1} Let $(M,L)$ be a Finsler manifold.
For a regular connection ${D}^{*}$ on $\pi^{-1}(TM)$ satisfying
conditions (a) and (b) of Theorem \ref{0th.h2},
 the (h)hv-torsion   of  $D^{*}$ coincides with
the (h)hv-torsion of the Cartan connection\,: $T^{*}=T$.
\end{lem}

\begin{proof}
The connection $D^{*}$ being regular, let  $h^{*}$ and $v^{*}$ be
the horizontal and vertical projectors associated with the
 decomposition (\ref{proj.}): $h^{*}=\beta^{*}\circ\rho, v^{*}=I-\beta^{*}\circ\rho
 $.\\
 As $D^{*}$ is non-metric with
 nonzero torsion ${\textbf{T}}^{*}$, then  ${D}^{*}$  is completely  determined by Equation (\ref{eq.r})
 with $D$ and $\textbf{T}$ replaced by $D^{*}$ and $\textbf{T}^{*}$ respectively.
\par
 Replacing $X, Y, Z$ by $\gamma \overline{X}, {h}^{*}Y, {h}^{*}Z $
 in (\ref{eq.r})
 and using conditions  (a) and  (b),  taking into account
  the fact that $\rho \circ \gamma=0$ and $\rho \circ h^{*}=\rho$, we get\vspace{-0.2cm}
\begin{equation}\label{eq.t11}
2g({D}^{*} _{\gamma \overline{X}}\rho Y,\rho Z) =\gamma
\overline{X}\cdot g(\rho Y,\rho Z)+
     g(\rho Y,\rho [h^{*}Z,\gamma \overline{X}])+g(\rho Z,\rho [\gamma \overline{X},h^{*}Y]).\vspace{-0.2cm}
\end{equation}

On the other hand, since the difference between two nonlinear
connections is a semi-basic vector form on $TM$  \cite{r21} and the
vertical distribution is completely integrable,
 then, noting that $\rho (V(TM))=0$, we get\vspace{-0.1cm}
\begin{equation}\label{04}
    \rho [h^{*}Z,\gamma \overline{X}]=\rho [hZ,\gamma
\overline{X}],\vspace{-0.1cm}
\end{equation}
 $h$ being the horizontal projector of the Cartan
connection (or the Barthel connection).
 \par
 Now, from (\ref{eq.t11}) and (\ref{04}), using  Theorem \ref{le.cc3}(i), we get\vspace{-0.2cm}
$${D}^{*} _{\gamma \overline{X}}\overline{Y}=\nabla _{\gamma
\overline{X}} \overline{Y}.\vspace{-0.1cm}$$
 Consequently, again by (\ref{04}), we get  ${T}^{*}( \overline{X},\overline{ Y})={T}(
\overline{X},\overline{ Y}).~$
\end{proof}

\noindent\textit{\textbf{Proof of theorem \ref{0th.h2}:}}
\par
By Lemma \ref{lem.ah1}, $ T^{*}(\overline{X}, \overline{\eta})=0$ as
$T$ has the same property. From which, together with  Proposition
\ref{eqv.}, it follows that $K^{*}=\gamma^{-1}$ on $V(TM)$ and the
associated nonlinear connection $\Gamma^{*}$ is given by
$\Gamma^{*}=\beta^{*}\circ\rho -\gamma \circ K^{*}$.
\par
Now, we prove that $\Gamma^{*}$ has the following properties:

 \vspace{4pt}
\noindent\textit{\textbf{$\Gamma^{*} $ is conservative:}} \
 $d_{h^{*}} E(X)=h^{*}X\cdot E=LD^{*}_{h^{*}X}L=0$, by condition
 (e).

 \noindent \textit{\textbf{$\Gamma^{*}$ is
homogenous :}} one can easily show that
  $$
    [\mathcal{C},v^{*}](X) =[\mathcal{C},\gamma(K^{*} X)]
    -\gamma (K^{*}[\mathcal{C},v^{*}X])-\gamma
    (K^{*}[\mathcal{C},h^{*}X]).
  $$
From which, using Lemma \ref{bracket2} and condition (d), noting
that $\gamma\,o\,K^{*}=id_{V(TM)},  K^{*}\circ \beta^{*}=0$ and
$[\mathcal{C},v^{*}X]$ is vertical, we get
\begin{eqnarray*}
    [\mathcal{C},v^{*}](X) &=& -\gamma (K^{*}[\mathcal{C},h^{*}X])=
    -\gamma \,o\,K^{*}([\gamma \overline{\eta},\beta^{*}\rho X])\\
    &=&-\gamma \,o\,K\{-
     \gamma(D^{*}_
     {\beta^{*} \rho X}\overline{\eta})
     +\beta^{*}( D^{*}_{\gamma \overline{\eta}}\rho X-T(\overline{\eta},\rho
     X))\}=0.
  \end{eqnarray*}
Therefore, $ [\mathcal{C},\Gamma^{*}]= -2[\mathcal{C},v^{*}]=0$ and
$\Gamma^{*}$ is thus  homogenous.

\vspace{4pt}
 \noindent
\textit{\textbf{$\Gamma^{*}$ is torsion-free:}} By the same argument
as in the proof of Theorem \ref{th.a}, we have
\begin{eqnarray*}
[J,v^{*}](X,Y)&=& J[\beta^{*} \rho X,\beta^{*} \rho Y]
   -v^{*}[\gamma \rho X,\beta^{*} \rho Y]+v^{*}[\gamma \rho Y,\beta^{*} \rho X].
\end{eqnarray*}
From which, together with Lemma \ref{bracket2} and condition (c), we
obtain
  \begin{eqnarray*}
  [J,v^{*}](X,Y)&=&J\{\gamma(\widehat{R^{*}}(\rho{X},\rho{Y}))
     + \beta^{*}(D^{*}_{h^{*}{X}}\rho{Y}-
     D^{*}_{h^{*}{Y}}\rho{X})\}\\
    &&-\gamma\circ K^{*}\{-
     \gamma(D^{*}_
     {h^{*}{Y}}\rho{X})
     +\beta^{*}( D^{*}_{J{X}}\rho{Y}-T(\rho{X},\rho{Y}))\}\\
     &&+\gamma\circ K^{*}\{-
     \gamma(D^{*}_
     {h^{*}{X}}\rho{Y})
     +\beta^{*}( D^{*}_{J{Y}}\rho{X}-T(\rho{Y},\rho{X}))\}=0.
     \end{eqnarray*}
    Hence
    $t:=\frac{1}{2}[J,\Gamma^{*}]=-[J,v^{*}]=0$.\\
   By the uniqueness of the Barthel connection
    (Theorem
  \ref{th.9a}), $\Gamma^{*} $ coincides with the Barthel
   connection $\Gamma=[J,G]$.\ \ $\Box$

\begin{thm} \label{th.h2} Let $(M,L)$ be a Finsler manifold and  $g$ the Finsler metric
defined by$\,L$. There exists a unique regular  connection ${D}^{*}$
on $\pi^{-1}(TM)$ such that
\begin{description}

  \item[(a)]  ${D}^{*}$ is vertically  metric\,{\em:} ${D}^{*}_{\gamma \overline{X}}\, g=0$,

   \item[(b)]  The (h)hv-torsion  $T^{*}$ of ${D}^{*}$ satisfies\,{\em:}
  $g({T}^{*}(\overline{X},\overline{Y}), \overline{Z})=
  g({T}^{*}(\overline{X},\overline{Z}),\overline{Y})$,

 \item[(c)]   The (h)h-torsion of ${D}^{*}$ vanishes\,{\em:} $Q^{*}=0
  $,

  \item[(d)]The (v)hv-torsion of ${D}^{*}$ vanishes\,{\em:}
  $\widehat{{P}^{*}}=0$,

\item[(e)] $D^{*}_{h^{*}X}L=0$.
  \end{description}
   \par Such a connection is called the Hashiguchi connection associated with the
Finsler manifold $(M,L)$.
\end{thm}

\begin{proof}
 First we prove the {\it \textbf{uniqueness}}.  The connection $D^{*}$ being regular,
 let $\Gamma^{*}$ be its associated nonlinear connection and
$h^{*}$ and $v^{*}$ its  horizontal and vertical projectors.\\
 By Theorem \ref{0th.h2}, we
 have $\Gamma^{*}=\Gamma=[J,G]$, and consequently
\begin{equation}\label{0000}
   v^{*}=v= \gamma \circ K , \ \ \ \  h^{*}=h=\beta \circ \rho,
      \ \ \ \ \ \ K^{*}=K, \ \ \ \ \beta^{*}=\beta.
\end{equation}
Also, by Lemma \ref{lem.ah1}, we have \ $T^{*}(K X,\rho Y)=T(K
X,\rho Y) $. \\
 Consequently,
\begin{equation}\label{qq1}
 {D}^{*}_{v{X}}\overline{Y}=\nabla _{v{X}}\overline{Y}.\vspace{-0.2cm}
\end{equation}

Now, using axiom (d), taking into account (\ref{0000}),  the
definition of ${{P}}^{*}$ and the identities $K\circ J=\rho$ and
$K\circ h=0$, we get \vspace{-0.2cm}
\begin{eqnarray*}
  0=\widehat{P}^{*}(\rho X,\rho Y)&=&\textbf{K}^{*}(hX,JY)\overline{\eta}
   =-{{D}}^{*}_{hX} {D}^{*}_{JY}\overline{\eta}
  +{{D}}^{*}_{JY} {D}^{*}_{hX}\overline{\eta} +{{D}}^{*}_{[hX,JY]}\overline{\eta} \\
  &=& - {{D}}^{*}_{hX}\rho Y+ K[hX,JY].\vspace{-0.3cm}
  \end{eqnarray*}
Hence,\vspace{-0.2cm}
 \begin{equation}\label{qq2}
  {{D}}^{*}_{hX}\rho Y= K[hX,JY]=\nabla_{[hX,JY]}\overline{\eta}=\nabla_
     {hX}\rho Y +\widehat{P}(\rho {X},\rho{Y}).
  \vspace{-0.2cm}
\end{equation}
Consequently,  (\ref{qq1}) and (\ref{qq2}) imply that
\vspace{-0.2cm}
\begin{equation}\label{qq3}
 {D}^{*}_{X}\overline{Y} = \nabla _{X}\overline{Y}
+ \widehat{P}(\rho X,\overline{Y}), \vspace{-0.2cm}
\end{equation}
which uniquely determines the connection $D^{*}$.
\par
\vspace{5pt}
 To prove the {\it \textbf{existence}} of $D^{*}$, we
define $D^{*}$ by (\ref{qq3}) and prove  the following properties\,:
\par
\vspace{7pt}
 \noindent \textit{\textbf{ ${D}^{*}$ satisfies
condition (a)}}\,: From  (\ref{qq3}), as  $\rho\circ\gamma=0$,
we get\\
 $ {\ \ \ }{D}^{*}_{\gamma \overline{X}}\overline{Y} =\nabla _{\gamma
\overline{X}}\overline{Y}$. Consequently,
 ${D}^{*}_{\gamma \overline{X}}\,g =\nabla
_{\gamma  \overline{X}}\,g=0,\ \forall \ \overline{X} \in \cp.$
\par
 \vspace{7pt}
 \noindent\textit{\textbf{ ${D}^{*}$ satisfies condition (b)}}\,:
 As  $D^{*}_{\gamma \overline{X}}\rho Y=\nabla_{\gamma \overline{X}}\rho
 Y$, taking into account (\ref{04}), one concludes that
${T}^{*}(\overline{X},\overline{Y})={T}(\overline{X},\overline{Y}).\
\ $
 $T^{*}$ has the property (b) as $T$ does.
\par
 \vspace{7pt}
 \noindent\textit{\textbf{ ${D}^{*}$ satisfies conditions (c)}}\,:
 Setting $\overline{Y}=\overline{\eta}$  into (\ref{qq3}), taking into account the fact that
  $\widehat{P}(\overline{X},\overline{\eta})=0$, it follows
that ${K}^{*}=K$. Since  $T^{*}(\overline{X},\overline{\eta})=0$ (as
${T}^{*}=T$), we have ${v}^{*}=\gamma\circ K^{*}=\gamma\circ K=v$
(by Proposition \ref{eqv.}). Consequently, $h^{*}=h$ and
$\beta^{*}=\beta$. Hence, from (\ref{qq3}) and the property that
$\widehat{P}$ is symmetric, one gets
\begin{eqnarray*}
 {Q}^{*}(\overline{X},\overline{ Y})&=&
 {D}^{*}_{\beta^{*}\overline{X}} \overline{Y}-{D}^{*}_{\beta^{*}\overline{Y}} \overline{X}
 -\rho[\beta^{*}\overline{X},\beta^{*}\overline{Y}]=
 {D}^{*}_{\beta\overline{X}} \overline{Y}-{D}^{*}_{\beta\overline{Y}} \overline{X}
 -\rho[\beta\overline{X},\beta\overline{Y}]\\
 &=&\nabla _{\beta\overline{X}}\overline{Y}
+ \widehat{P}(\overline{X},{\overline{Y}})-\nabla
_{\beta\overline{Y}}{\overline{X}} -
\widehat{P}(\overline{Y},{\overline{X}})-\rho[\beta\overline{X},\beta\overline{Y}]=
Q(\overline{X},\overline{Y})=0\vspace{-0.2cm}.
\end{eqnarray*}
\par
 \noindent\textit{\textbf{ ${D}^{*}$ satisfies condition (d)}}\,: As
 $K\circ \beta=0$ and $\widehat{P}(\overline{X},\overline{\eta})=0$, using (\ref{qq3}), we get\vspace{-0.2cm}
\begin{eqnarray*}
  \widehat{{P}^{*}}(\overline{X},\overline{Y})
  &=&{{P}}^{*}(\overline{X},\overline{Y})\overline{\eta}=
  \textbf{{{K}}}^{*}(\beta\overline{X},\gamma\overline{Y})\overline{\eta}\\
  &=&-{{D}}^{*}_{\beta\overline{X}}D^{*}_{\gamma\overline{Y}}\overline{\eta}
  +{{D}}^{*}_{\gamma\overline{Y}}D^{*}_{\beta\overline{X}}\overline{\eta}
+{{D}}^{*}_{[\beta\overline{X},\gamma\overline{Y}]}\overline{\eta}=-{{D}}^{*}_{\beta\overline{X}}
\nabla_{\gamma\overline{Y}}\overline{\eta}+
\nabla_{[\beta\overline{X},\gamma\overline{Y}]}\overline{\eta}\\
&=&-\nabla_{\beta\overline{X}}
\nabla_{\gamma\overline{Y}}\overline{\eta}-
\widehat{P}(\overline{X},\nabla_{\gamma\overline{Y}}\overline{\eta})
+\nabla_{[\beta\overline{X},\gamma\overline{Y}]}\overline{\eta}\\
&=&\{-\nabla_{\beta\overline{X}}
\nabla_{\gamma\overline{Y}}\overline{\eta}
+\nabla_{[\beta\overline{X},\gamma\overline{Y}]}\overline{\eta}\}-
\widehat{P}(\overline{X},\overline{Y})=0.
\end{eqnarray*}
 \par
\vspace{5pt}
 \noindent \textit{\textbf{ ${D}^{*}$ satisfies condition(e)}}\,:
As $h^{*}=h$, then    $$ L\,{D}^{*}_{{h^{*}}X}L=
{D}^{*}_{{h}X}E=hX\cdot E=2g(\nabla_{hX}\overline{\eta},\overline{
\eta})=0.$$
 This complete the proof.
\end{proof}

It is to be  noted  that  when we localize the above result, the {
\it{local expressions}} obtained coincide with the classical
expressions found in \cite{r34}, \cite{r23}, \cite{r93}...etc. (c.f.
Appendix 2).
\par
\vspace{5pt}
In view of the above theorem taking Theorem \ref{th.5}
into account, we have\,:
 \vspace{-0.2cm}
\begin{thm}\label{th.t12} The Hashiguchi  connection $ {D}^{*} $  is
given in terms of the Cartan connection \emph{(}{\it or the Berwald
connection}\emph{)} by\,\emph{:} \vspace{-0.2cm}
  \begin{eqnarray}\label{4}
   {D}^{*}_{X}\overline{Y} = \nabla _{X}\overline{Y}
+ {\widehat{P}}(\rho X,\overline{Y}) ={D}^{\circ}_{X}\overline{Y}
+{T}(K X, \overline{Y}) . \vspace{-0.2cm}
  \end{eqnarray}
In particular, we have
\begin{description}
  \item[(a)] $ {D}^{*}_{\gamma \overline{X}}\overline{Y}=\nabla _{\gamma
  \overline{X}}\overline{Y}={D}^{\circ}_{\gamma \overline{X}}\overline{Y}
+{T}( \overline{X}, \overline{Y})$.

 \item[(b)] $ {D}^{*}_{\beta \overline{X}}\overline{Y}=\nabla _{\beta
  \overline{X}}\overline{Y}+\widehat{P}(\overline{X},\overline{Y})
  ={D}^{\circ}_{\beta \overline{X}}\overline{Y}.$
\end{description}
\end{thm}

\begin{cor}\label{pp.t6}For every $\overline{X},
\overline{Y}, \overline{Z}\in\cp$, we have
\begin{description}
  \item[(a)]$ {D}^{*}_{\gamma \overline{X}}g=0$,
  \item[(b)]$ ({D}^{*}_{\beta
  \overline{X}}g)(\overline{Y},\overline{Z})=-2g(\widehat{P}(\overline{X},
  \overline{Y}),\overline{Z})$.
  \end{description}
  Consequently, $(M,L)$ is a Landsberg manifold  if, and only if,
 ${D}^{*}_{\beta \overline{X}}g=0 $.
\end{cor}

Concerning the curvature tensors of the Hashiguchi connection
${D}^{*}$, we have\,:\vspace{-0.2cm}
\begin{prop}\label{pp.tt3} The curvature tensor ${\textbf{K}}^{*}$ of the Hashiguchi connection
${D}^{*}$ is given in terms of the  curvature tensor ${\textbf{K}}$
of the Cartan connection $\nabla$ by\,:\vspace{-0.1cm}
\begin{equation*}
\left.
    \begin{array}{rcl}
    {\textbf{K}}^{*}(X,Y)\overline{Z} &=&{\textbf{K}}(X,Y)\overline{Z}-\widehat{P}({\textbf{T}}(X,Y), \overline{Z})
   -\mathfrak{U}_{X,Y}\{(\nabla_{X}\widehat{P})(\rho Y,\overline{Z})
   +\widehat{P}(\rho X,\widehat{P}(\rho Y,\overline{Z}))\}.\vspace{-0.1cm}
\end{array}
  \right.
\end{equation*}
In particular, we have\vspace{-0.1cm}
\begin{description}

  \item[(a)]${S}^{*}(\overline{X},\overline{Y})\overline{Z}= S(\overline{X},\overline{Y})\overline{Z}$.

  \item[(b)]${P}^{*}(\overline{X},\overline{Y})\overline{Z}= P(\overline{X},\overline{Y})\overline{Z}
  + \widehat{P}({{T}}(\overline{X},\overline{Y}), \overline{Z})+(\nabla_{\gamma \overline{Y}}\widehat{P})(\overline{X},\overline{Z})$.

  \item[(c)]${R}^{*}(\overline{X},\overline{Y})\overline{Z}= R(\overline{X},\overline{Y})\overline{Z}-
  \mathfrak{U}_{\overline{X},\overline{Y}}\{ (\nabla_{\beta \overline{X}}\widehat{P})(\overline{Y}, \overline{Z})
  +\widehat{P}(\overline{X},\widehat{P}(\overline{Y},\overline{Z})) \}$.
\end{description}
\end{prop}
\vspace{0.2cm}

\begin{cor}~\par
\begin{description}
     \item[(a)] $\widehat{{S}^{*}}(\overline{X},\overline{Y})= 0$.
    \item[(b)]$\widehat{P^{*}}(\overline{X},\overline{Y})=0$,\,\,
    ${P}^{*}(\overline{X},\overline{\eta})\overline{Y}=0$,\,\,
     ${P}^{*}(\overline{\eta},
     \overline{X})\overline{Y}=-\widehat{P}(\overline{X},\overline{Y})$.
    \item[(c)] $ \widehat{{R}^{*}}(\overline{X},\overline{Y})= \widehat{R}(\overline{X},\overline{Y})$.
\end{description}
\end{cor}

We terminate the paper by some comments\,:
\begin{itemize}
\item  \textit{On a Finsler manifold $(M,L)$, there are
\it{canonically} associated  four linear connections\,: the Cartan
connection $\nabla$ (Thm. \ref{le.cc3}), the Chern connection
$D^{\diamond}$ (Thm. \ref{th.r1}), the Hashiguchi connection $ D^{*}
$ (Thm. \ref{th.h2}) and   the  Berwald connection $D^{\circ} $
(Thm. \ref{bth2.h2}). The nonlinear connection associated with each
of these linear connections coincides with  the Barthel connection.}

\item  \textit{There are two
 methods by which a given Finsler connection is converted to some other connection.
 From Theorem
 \ref{th.r5}, we observe that  the Chern connection $D^{\diamond}$ is obtained
  from the Cartan connection $\nabla$ by subtracting the
  (h)hv-torsion  $T$ from  its vertical counterpart. Moreover,
  from Theorem \ref{th.t12}, we observe that  the Hashiguchi
   connection $D^{*}$ is obtained from the Cartan connection
   by adding the (v)hv-torsion $\widehat{P}$ to
  its horizontal counterpart.
The former process is called the $C$-process  and the latter process
is called the $P^{1}$-process by M. Matsumoto\,\cite{r34}}.

\textit{Interestingly, the two processes commute with each other. If
we apply to the Cartan connection   the $C$-process followed by the
$P^{1}$-process, we obtain the Berwald connection $D^{\circ}$
(passing by the Chern connection) (cf. Theorem \ref{th.5}). On the
other hand, if we apply the $P^{1}$-process followed by the
$C$-process, we also obtain the Berwald connection (passing by the
Hashiguchi connection). That is},
$$ D^{\circ} \stackrel{C\text{-process}}
\longleftarrow D^{*} \stackrel{P^{1}\text{-process}} \longleftarrow
\nabla \stackrel{C\text{-process}} \longrightarrow
D^{\diamond}\stackrel{P^{1}\text{-process}}\longrightarrow D^{\circ}
$$
\end{itemize}

\vspace{9pt}
\newpage
\begin{center}
{\bf{\Large{ Appendix 1. Intrinsic Comparison}}}
\end{center}
\par The following tables establish a concise comparison concerning the
canonical linear connections in Finsler geometry as well as the
fundamental geometric objects associated with them.
 \begin{center}{\bf{Table 1.}}
\end{center}
 {\tiny
\begin{center}
\begin{tabular}
{|c|c|c|c|c|}\hline
&&&&\\
 {\bf connection} &{\bf Cartan:\,$\nabla$ }& {\bf Chern:\,$D^{\diamond}$ } &{\bf Hashiguchi:\,$D^{*}$ }
 &{\bf Berwald:\,$D^{\circ}$}
\\[0.1 cm]\hline
&&&&\\
{\bf v-counterpart} & $\nabla _{\gamma\overline{X}}\overline{Y}$&
 $D^{\diamond} _{\gamma\overline{X}}
      \overline{Y}=\nabla _{\gamma\overline{X}}
      \overline{Y}-T(\overline{X},\overline{Y})$&
${D}^{*} _{\gamma\overline{X}}
      \overline{Y}=\nabla _{\gamma\overline{X}}
      \overline{Y}$&
${D}^{\circ} _{\gamma\overline{X}}
      \overline{Y}=\nabla _{\gamma\overline{X}}
      \overline{Y}-T(\overline{X},\overline{Y})$
\\[0.1 cm]
{\bf h-counterpart} & $\nabla _{\beta\overline{X}}\overline{Y}$&
$D^{\diamond} _{\beta\overline{X}}
      \overline{Y}=\nabla _{\beta\overline{X}}
      \overline{Y}$&
${D}^{*} _{\beta\overline{X}}
      \overline{Y}= \nabla _{\beta\overline{X}}
      \overline{Y}+\widehat{P}(\overline{X},\overline{Y})$&
      ${D}^{\circ} _{\beta\overline{X}}
      \overline{Y}=\nabla _{\beta\overline{X}}
      \overline{Y}+\widehat{P}(\overline{X},\overline{Y})$

\\[0.1 cm]\hline
&&&&\\
{\bf (h)v-torsion} & $0$& $0$& $0$&$0$
\\[0.1 cm]
{\bf (h)hv-torsion} & $T$& $0$& $T$&$0$

\\[0.1 cm]
{\bf (h)h-torsion} & $0$& $0$&$0$&$0$
\\[0.1 cm]\hline
&&&&\\
{\bf (v)v-torsion} & $0$& $0$& $0$&$0$
\\[0.1 cm]
{\bf (v)hv-torsion}& $\widehat{P}=\nabla_{G}T$&
$\widehat{P}$&$0$&$0$
\\[0.1 cm]
{\bf (v)h-torsion} & ${\widehat{R}}=-K\mathfrak{R}$&
${\widehat{R}}$& ${\widehat{R}}$&${\widehat{R}}$
\\[0.1 cm]\hline
&&&&\\
{\bf v-curvature} & ${S}$& $0$& $S$&$0$
\\[0.1 cm]
{\bf hv-curvature} & ${P}$& ${P^{\diamond}}$&
${P^{*}}$&${P^{\circ}}$
\\[0.1 cm]
{\bf h-curvature} & ${R}$& $R^{\diamond}$& ${R^{*}}$&${R^{\circ}}$
\\[0.1 cm]\hline
&&&&\\
{\bf v-metricity}& $\nabla _{\gamma\overline{X}}{g}=0$&
$D^{\diamond} _{\gamma\overline{X}}{g}=2g(T(\overline{X},.),.)$&
$D^{*} _{\gamma\overline{X}}{g}=0$&$D^{\circ}
_{\gamma\overline{X}}{g}=2g(T(\overline{X},.),.)$
\\[0.1 cm]
{\bf h-metricity}& $\nabla _{\beta\overline{X}}{g}=0$& $D^{\diamond}
_{\beta\overline{X}}{g}=0$& $D^{*}
_{\beta\overline{X}}{g}=-2g(\widehat{P}(\overline{X},.),.)$&$D^{\circ}
_{\beta\overline{X}}{g}=-2g(\widehat{P}(\overline{X},.),.)$
\\[0.1 cm]\hline
\end{tabular}
\end{center}
}

\begin{center}{\bf{Table 2.}}
\\[0.2cm]
\begin{tabular}
{|c|c|}\hline
&\\
{\bf connection}&\bf{curvature tensors}
\\[0.1 cm]\hline
&\\
{\bf Cartan }&  v-curvature\,:\,
$S(\overline{X},\overline{Y})\overline{Z}:=-\nabla_{\gamma
\overline{X}} \nabla_{\gamma\overline{Y}}
\overline{Z}+\nabla_{\gamma \overline{Y}}
\nabla_{\gamma\overline{X}}
\overline{Z}+\nabla_{[\gamma\overline{X},\gamma\overline{X}]}
\overline{Z}$.
\\
\ \ & hv-curvature\,:\,
$P(\overline{X},\overline{Y})\overline{Z}:=-\nabla_{\beta
\overline{X}} \nabla_{\gamma\overline{Y}}
\overline{Z}+\nabla_{\gamma \overline{Y}} \nabla_{\beta\overline{X}
} \overline{Z}+\nabla_{[\beta\overline{X},\gamma\overline{X}]}
\overline{Z}$.
\\
\ \ & h-curvature\,:\,
$R(\overline{X},\overline{Y})\overline{Z}:=-\nabla_{\beta
\overline{X}} \nabla_{\beta\overline{Y}} \overline{Z}+\nabla_{\beta
\overline{Y}} \nabla_{\beta\overline{X} }
\overline{Z}+\nabla_{\beta\overline{X},\beta\overline{X}]}
\overline{Z}$.
\\[0.1 cm]\hline
&\\
{\bf Chern  }& $S^{\diamond}(\overline{X},\overline{Y})\overline{Z}=
0$.
\\
\ \ & $P^{\diamond}(\overline{X},\overline{Y})\overline{Z}=
  P(\overline{X},\overline{Y})\overline{Z}-T(\widehat{P}(\overline{X},\overline{Y}), \overline{Z})
  +(\nabla_{\beta \overline{X}}T)(\overline{Y},\overline{Z})$.
\\
\ \ &$R^{\diamond}(\overline{X},\overline{Y})\overline{Z}=
  R(\overline{X},\overline{Y})\overline{Z}-T(\widehat{R}(\overline{X},\overline{Y}),
  \overline{Z})$.
\\[0.1 cm]\hline
&\\
{\bf Hashiguchi }& ${S}^{*}(\overline{X},\overline{Y})\overline{Z}=
S(\overline{X},\overline{Y})\overline{Z}$.
\\
\ \ & ${P}^{*}(\overline{X},\overline{Y})\overline{Z}=
P(\overline{X},\overline{Y})\overline{Z}
  + \widehat{P}({{T}}(\overline{X},\overline{Y})
  , \overline{Z})+(\nabla_{\gamma \overline{Y}}\widehat{P})(\overline{X},\overline{Z})$.
\\
\ \ & ${R}^{*}(\overline{X},\overline{Y})\overline{Z}=
R(\overline{X},\overline{Y})\overline{Z}-
  \mathfrak{U}_{\overline{X},\overline{Y}}\{
   (\nabla_{\beta \overline{X}}\widehat{P})(\overline{Y}, \overline{Z})
  +\widehat{P}(\overline{X},\widehat{P}(\overline{Y},\overline{Z})) \}$.
\\[0.1 cm]\hline
&\\
{\bf Berwald}&
${{S}}^{\circ}(\overline{X},\overline{Y})\overline{Z}= 0 $.
\\
\ \ &${{P}}^{\circ}(\overline{X},\overline{Y})\overline{Z}=
P(\overline{X},\overline{Y})\overline{Z}+
  (\nabla_{\gamma \overline{Y}}\widehat{P})(\overline{X},\overline{Z})
  +\widehat{P}(T(\overline{Y},\overline{X}),\overline{Z})+$
  \\
\ \ &$ { \ \  \ \ \ }
   +\widehat{P}(\overline{X},T(\overline{Y},\overline{Z}))
   +(\nabla_{\beta\overline{X}}T)(\overline{Y},\overline{Z})-$
\\
\ \ & ${ \ \ \ } -T(\overline{Y},
\widehat{P}(\overline{X},\overline{Z}))-T( \widehat{P}(\overline{X},
\overline{Y}), \overline{Z})$.
\\
\ \ & ${{R}}^{\circ}(\overline{X},\overline{Y})\overline{Z}=
R(\overline{X},\overline{Y})\overline{Z}-
T(\widehat{R}(\overline{X},\overline{Y}),\overline{Z})-
\mathfrak{U}_{\overline{X},\overline{Y}}\{(\nabla_{\beta
\overline{X}}\widehat{P})(\overline{Y}, \overline{Z})$
\\
\ \ &${ \qquad  \ \ \ \ \ \ }
+\widehat{P}(\overline{X},\widehat{P}(\overline{Y},\overline{Z}))\}.$
\\[0.1 cm]\hline
\end{tabular}
\end{center}


\bigskip

\bigskip
\begin{center}
{\bf\Large{\textbf{Appendix 2. Local Comparison}}}\end{center}

 \vspace{3pt}
 \par For the sake of completeness, we present in this appendix
 a brief and concise survey of the local expressions of the most
 important geometric objects treated in the paper.
\vspace{5pt}
\par
  Let $(U,(x^{i}))$ be  a system  of local coordinates on
 $M$ and $(\pi^{-1}(U),(x^i,y^i))$ the associated system of local coordinates on $TM$.
 We use the following notations\,:\\
 $(\pa_{i}):=(\frac{\pa}{\pa x^i})$: the natural basis of $T_{x}M,\, x\in
 M$,\\
 $(\paa_{i}):=(\frac{\pa}{\pa y^i})$: the natural basis of $V_{u}(\T M),\, u\in
 \T M$,\\
$(\pa_{i},\paa_{i})$: the natural basis of $T_{u}(\T M)$,\\
$(\overline{\pa}_{i} )$: the natural basis of the fiber over $u$ in
$\p$ ($\overline{\pa}_{i} $ is the lift of $\pa_{i}$ at $u$).
\vspace{4pt}
\par To  a Finsler manifold $(M,L)$, we associate the
geometric
objects\,:\\
$g_{ij}:= \frac{1}{2} \paa_{i} \paa_{j}L^{2}= \paa_{i} \paa_{j}E$:
the
Finsler metric tensor,\\
$G^h$: the components of the canonical spray,\\
$G^{h}_{i}:=\paa_{i}\,G^h$,\\
$G^{h}_{ij}:=\paa_{j}\,G^h_{i}=\paa_{j}\paa_{i}\,G^h$,\\
$(\delta_{i}):=(\pa_{i}-G^{h}_{i}\,\paa_{h})$: the basis of
$H_{u}(\T
M)$ adapted to $G^{h}_{i}$,\\
$(\delta_{i}, \paa_{i})$: the basis of $T_{u}(\T M)=H_{u}(\T
M)\oplus V_{u}(\T M)$ adapted to $G^{h}_{i}$.
\vspace{5pt}
\par
We have\,:\\
$\gamma(\overline{\pa}_{i})=\paa_{i}$,\\
$\rho(\pa_{i})=\overline{\pa}_{i}$, \ \ $\rho(\paa_{i})=0$,\ \
$\rho(\delta_{i})=\overline{\pa}_{i}$,\\
$\beta(\overline{\pa}_{i})=\delta_{i}$,\\
$J(\pa_{i})= \paa_{i}$, \ \ $J(\paa_{i})=0$,\ \
$J(\delta_{i})= \paa_{i}$,\\
$h:= \beta o \rho= dx^{i} \otimes \pa_{i}- G^{i}_{j}\, dx^{j}
\otimes \paa_{i}$,\ \ \ \   $v:=\gamma o K=dy^{i} \otimes \paa_{i}+
G^{i}_{j}\, dx^{j} \otimes \paa_{i} $.
\vspace{5pt}
\par
 We define\,:\\
\noindent $\gamma^{h}_{ij}:= \frac{1}{2}\,g^{h\ell}(\pa_{i}\,g_{\ell
j}+\pa_{j}\,g_{i\ell }- \pa_{\ell}\,g_{i j}
 )$,\\
  $C^{h}_{ij}:= \frac{1}{2}\,g^{h\ell}(\paa_{i}\,g_{\ell j}+\paa_{j}\,g_{i\ell }-
   \paa_{\ell}\,g_{i j})=
  \frac{1}{2}\,g^{h\ell}\,\paa_{i}\,g_{\ell j}$,\\
  $\Gamma^{h}_{ij}:= \frac{1}{2}\,g^{h\ell}(\delta_{i}\,g_{\ell j}+\delta_{j}\,g_{i\ell }-
  \delta_{\ell}\,g_{i j})$.
 \vspace{5pt}
\par
 For a Finsler connection  $F\Gamma$:
  $(\ \textbf{F}^{h}_{ij}, \textbf{N}^{h}_{i}, \ \textbf{C}^{h}_{ij})$, we have two covariant derivatives\,:\\
    The $h$-covariant derivative $A^{i}_{j|k}$  given by
    $A^{i}_{j|k}:=\delta_{k}\,A^{i}_{j}+A^{m}_{j}\,\textbf{F}^{i}_{mk}
   -A^{i}_{m}\,\textbf{F}^{m}_{jk}$,  \\
 The  $v$-covariant derivative $A^{i}_{j}|_{k}$ given by
    $A^{i}_{j}|_{k}:=\paa_{k}\,A^{i}_{j}+A^{m}_{j}\,\textbf{C}^{i}_{mk}
   -A^{i}_{m}\,\textbf{C}^{m}_{jk}$.
  \vspace{5pt}
\par
The  canonical spray and the canonical connections in Finsler geometry are as follows\,:\\
\vspace{5pt}
 -- \ The canonical spray $G$:
$G^h=\frac{1}{2}\,\gamma^{h}_{ij}\,\,y^i
   \,y^j$,\\
\vspace{4pt}
-- \ The Barthel connection  $\Gamma$: $G^{h}_{i}= \paa_{i}\,G^{h}$,\\
\vspace{4pt}
-- \  The Cartan connection  $C\Gamma$:
  $(\ \Gamma^{h}_{ij}, G^{h}_{i}, C^{h}_{ij})$,\\
\vspace{4pt}
 -- \  The Chern (Rund) connection  $R\Gamma$:
  $( \ \Gamma^{h}_{ij}, G^{h}_{i},0)$, \\
\vspace{4pt}
-- \   The Hashiguchi connection  $H\Gamma$: $(\ G^{h}_{ij}, G^{h}_{i}, C^{h}_{ij})$,\\
\vspace{4pt}
 -- \   The Berwald connection  $B\Gamma$: $(\
G^{h}_{ij}, G^{h}_{i},\ 0)$, \\
where $G^{h}_{ij}$ is given by\,:  $G^{h}_{ij}=
\paa_{j}\,\paa_{i}\,G^{h}= \Gamma^{h}_{ij}+ C^{h}_{ij}{_{|k}}\,y^k=
\Gamma^{h}_{ij}+ C^{h}_{ij}{_{|0}}\,
;\,\,C^{h}_{ij}{_{|0}}:=C^{h}_{ij}{_{|k}}y^k $.

\vspace{7pt}
\par
In the next table, we give the local expressions for the fundamental
tensors associated with the above linear connections.

\begin{center}{\bf{Table 3: Summary of local expressions}}
\\[0.4 cm]

\begin{tabular}
{|c|c|c|c|c|}\hline
&&&&\\
 \ {\bf connection } &{\bf Cartan }& {\bf Chern  } &{\bf Hashiguchi }&{\bf Berwald}
\\[0.1 cm]\hline
&&&&\\
$(\textbf{F}^{h}_{ij}, \textbf{N}^{h}_{i},\textbf{C}^{h}_{ij})$&
$(\Gamma^{h}_{ij}, G^{h}_{i}, C^{h}_{ij})$& $( \Gamma^{h}_{ij},
G^{h}_{i}, 0)$& $( G^{h}_{ij}, G^{h}_{i}, C^{h}_{ij})$&$(
G^{h}_{ij}, G^{h}_{i}, 0)$
\\[0.1 cm]\hline
&&&&\\
{\bf (h)v-torsion} & $0$& $0$& $0$&$0$
\\[0.1 cm]
{\bf (h)hv-torsion} & $C^{i}_{jk}$& $0$& $C^{i}_{jk}$&$0$

\\[0.1 cm]
{\bf (h)h-torsion} & $0$& $0$&$0$&$0$
\\[0.1 cm]\hline
&&&&\\
{\bf (v)v-torsion} & $0$& $0$& $0$&$0$
\\[0.1 cm]
{\bf (v)hv-torsion} &${P}^{i}_{jk}=C^{i}_{jk|0}$&
$P^{i}_{jk}$&$0$&$0$
\\[0.1 cm]
{\bf (v)h-torsion} &
${R}^{i}_{jk}=\delta_{k}\,G^{i}_{j}-\delta_{j}\,G^{i}_{k}$&
${R}^{i}_{jk}$& $R^{i}_{jk}$&${R}^{i}_{jk}$
\\[0.1 cm]\hline
&&&&\\
{\bf v-curvature} & ${S}^{h}_{ijk}$& $0$& $S^{h}_{ijk}$&$0$
\\[0.1 cm]
{\bf hv-curvature} & ${P}^{h}_{ijk}$& ${P}^{\diamond h}_{ijk}$&
${P}^{* h}_{ijk}$&${P}^{\circ h}_{ijk}$
\\[0.1 cm]
{\bf h-curvature} & ${R}^{h}_{ijk}$& ${R}^{\diamond h}_{ijk}$&
${R}^{* h}_{ijk}$&${R}^{\circ h}_{ijk}$
\\[0.1 cm]\hline
&&&&\\
{\bf v-metricity}& ${g}_{ij}|_{k}=0$&
${g}_{ij}\hspace{-0.1cm}\stackrel{\diamond}|_{k}=2C_{ijk}$&
${g}_{ij}\hspace{-0.1cm}\stackrel{*}|_{k}=0$&${g}_{ij}\hspace{-0.1cm}\stackrel{\circ}|_{k}=2C_{ijk}$
\\[0.1 cm]
{\bf h-metricity}& ${g}_{ij|k}=0$&
${g}_{ij\stackrel{\diamond}|k}=0$&
${g}_{ij\stackrel{*}|k}=-2P_{ijk}$&${g}_{ij\stackrel{\circ}|k}=-2P_{ijk}$
\\[0.1 cm]\hline
&&&&\\
{\bf v-cov. derivative}& ${A}^{i}_{j}|_{k}$&
${A}^{i}_{j}\hspace{-0.1cm}\stackrel{\diamond}|_{k}=\paa_{k}\,{A}^{i}_{j}$&
${A}^{i}_{j}\hspace{-0.1cm}\stackrel{*}|_{k}={A}^{i}_{j}|_{k}$&${A}^{i}_{j}\hspace{-0.1cm}\stackrel{\circ}|_{k}=\paa_{k}\,{A}^{i}_{j}$
\\[0.1 cm]
{\bf h-cov. derivative}& ${A}^{i}_{j|k}$&
${A}_{j\stackrel{\diamond}|k}^{i}={A}_{j|k}^{i}$&
${A}^{i}_{j\stackrel{*}|k}$&${A}^{i}_{j\stackrel{\circ}|k}={A}^{i}_{j\stackrel{*}|k}$
\\[0.1 cm]\hline
\end{tabular}
\end{center}
\vspace{7pt}
The curvature tensors appeared in the above table are given by\,:\vspace{9pt}\\
 \noindent$R^{i}_{hjk}=K^{i}_{hjk}+C^{i}_{hm}
R^{m}_{jk}$,\quad
$P^{i}_{hjk}=\paa_{k}\Gamma^{i}_{hj}-C^{i}_{hk|j}+C^{i}_{hm}\,P^{m}_{jk}$,\quad
$S^{i}_{hjk}=\mathfrak{U}_{jk}\{ C^{m}_{hk}\,C^{i}_{mj}\}$,\vspace{7pt}\\
$R^{\diamond i}_{hjk}=K^{i}_{hjk}$,\quad $P^{\diamond
i}_{hjk}=\paa_{k}\,\Gamma^{i}_{hj}$,\vspace{7pt}\\
$R^{*i}_{hjk}=\mathfrak{U}_{jk}\{
\delta_{k}\,G^{i}_{hj}+G^{m}_{hj}\, G^{i}_{mk}\}+C^{i}_{hm}\,
R^{m}_{jk}$,\quad
 $P^{*i}_{hjk}=\paa_{k}\,G^{i}_{hj}-C^{i}_{hk\stackrel{*}|j}$,\quad
 $S^{*i}_{hjk}=S^{i}_{hjk}$,\vspace{2pt}\\
 ${{R}}^{\circ i}_{hjk}=\mathfrak{U}_{jk}\{
\delta_{k}\,G^{i}_{hj}+G^{m}_{hj}\, G^{i}_{mk}\}$,\quad
${{P}}^{\circ i}_{hjk}=\paa_{k}\,G^{i}_{hj}=:G^{i}_{hjk}$,\vspace{9pt}\\
 where $K^{i}_{hjk}=\mathfrak{U}_{jk}\{
\delta_{k}\,\Gamma^{i}_{hj}+\Gamma^{m}_{hj} \,\Gamma^{i}_{mk}\}\,$
and $\,\mathfrak{U}_{jk}A_{jk}:=A_{jk}-A_{kj}$.

\bigskip

\bigskip
\providecommand{\bysame}{\leavevmode\hbox
to3em{\hrulefill}\thinspace}
\providecommand{\MR}{\relax\ifhmode\unskip\space\fi MR }
\providecommand{\MRhref}[2]{%
  \href{http://www.ams.org/mathscinet-getitem?mr=#1}{#2}
} \providecommand{\href}[2]{#2}

\end{document}